\documentclass[11pt,a4paper]{amsart}
\usepackage{amsmath,amsfonts,amssymb,epsf,mathrsfs}
\usepackage[all]{xy}

\usepackage[T1]{fontenc}
\usepackage[backref=page,bookmarksopen=true]{hyperref}

\let\bbbibitem\bibitem
\renewcommand{\bibitem}[2][]{\bbbibitem[#1]{#2}\label{#2}}

\def\fin{\hfill\hbox{\hskip .2cm $\square$}\medskip}

\newtheorem{theo}{Theorem}[section]
\newtheorem{lemma}[theo]{Lemma}
\newtheorem{prop}[theo]{Proposition}
\newtheorem{coro}[theo]{Corollary}
\newtheorem{conj}[theo]{Conjecture}
\newtheorem{defi}[theo]{Definition}

\newenvironment{demo}[1][\hspace{-3pt}]{{\noindent\em Proof #1.~ }}{\fin}

\theoremstyle{remark}
\newtheorem{rema}[theo]{Remark}

\def\cal{\mathcal}
\def\a{\alpha}
\def\b{\beta}

\def\g{\gamma}
\def\G{\Gamma}
\def\Th{\Theta}
\def\R{{\mathbb R}} 

\def\Ad{{\rm Ad}}
\def\ad{{\rm ad}}

\def\PU{{\rm PU}}

\def\GL{{\rm GL}}
\def\PGL{{\rm PGL}}

\def\Hom{{\rm Hom}}
\def\U{{\rm U}}

\def\SO{{\rm SO}}

\def\Ker{{\rm Ker\,}}

\def\C{{\mathbb C}}

\def\N{{\mathbb N}}
\def\Z{{\mathbb Z}}
\def\R{{\mathbb R}}
\def\fd{\longrightarrow}

\def\la{\langle}
\def\ra{\rangle}

\def\om{\omega}

\def\rk{{\rm rk}\,}

\def\t{\theta}

\def\o{\omega}

\def\U{{\rm U}}

\def\gg{{\mathfrak g}}
\def\kg{{\mathfrak k}}

\def\pg{{\mathfrak p}}
\def\hg{{\mathfrak h}}
\def\qg{{\mathfrak q}}
\def\ug{{\mathfrak u}}

\def\ng{{\mathfrak n}}
\def\vg{{\mathfrak v}}
\def\sg{{\mathfrak s}}
\def\lg{{\mathfrak l}}

\def\P{{\mathbb P}}

\def\trans{{}^t\hspace{-2pt}}

\def\E{{\mathbb E}}

\newcommand{\lo}{\longrightarrow}
\newcommand{\sspace}{\vspace{0.1cm}}

\newcommand{\z}{{\mathbb Z}}
 
\newcommand{\nn}{{\mathbb N}} 
\newcommand{\real}{{\mathbb R}}
\newcommand{\cl}{{\mathbb C}}

\newcommand{\proj}{\mathbf P}

\newcommand{\Ga}{\Gamma}

\newcommand{\RR}{{R}}

\newcommand{\GG}{{G}}
\newcommand{\UU}{{U}}
\newcommand{\A}{{A}}

\newcommand{\hooklongrightarrow}{\lhook\joinrel\longrightarrow}
\newcommand{\twoheadlongrightarrow}{\relbar\joinrel\twoheadrightarrow}

\begin{document}

\title[On the second cohomology of K\"ahler groups]{On the second cohomology of K\"ahler groups}

\author{Bruno Klingler} 
\author{Vincent Koziarz} 
\author{Julien Maubon}
\address{IMJ, Universit\'e Paris-Diderot - Paris VII, CNRS, 175 rue du chevaleret, F-75013 Paris, France}
\email{klingler@math.jussieu.fr}
\address{IECN, Nancy-Universit\'e, CNRS, INRIA, Boulevard des Aiguillettes B. P. 239, F-54506 
Vand\oe uvre-l\`es-Nancy, France}
\email{koziarz@iecn.u-nancy.fr}  
\email{maubon@iecn.u-nancy.fr}

\date{\today}

\sloppy

\begin{abstract}  
Carlson and Toledo conjectured that if an infinite group $\G$ is the fundamental group of a
compact K\"ahler manifold, then virtually $H^2(\G, \R)\not =0$. We assume that $\G$ admits an unbounded reductive rigid linear representation. This representation necessarily comes from a complex variation of Hodge structure ($\C$-VHS) on the K\"ahler manifold. We prove the conjecture under some assumption on the $\C$-VHS. We also study some related geometric/topological properties of period domains associated to such $\C$-VHS.   
\end{abstract}

\maketitle
\section{Introduction} \label{intro}

The general setting for this paper is the study of topological
properties of compact K\"ahler manifolds, in particular complex smooth
projective varieties. The possible homotopy types for these spaces are
essentially unknown (cf. \cite{si4}).
Here we restrict ourselves to an a priori simpler
question: the study of fundamental groups of compact K\"ahler
manifolds (the so-called K\"ahler groups).

It is well-known that any finitely presented group $\Ga$ can be realized as
the fundamental group of a $4$-dimensional compact real manifold, or
even of a symplectic complex compact $3$-fold. A classical result of Serre shows
that any finite group can be realized as a projective group (i.e. the fundamental group of a complex smooth
projective variety). On the other hand there are many known
obstructions for an infinite finitely presented group being K\"ahler (we refer to \cite{abckt}
for a panorama). Most of them come from Hodge theory (Abelian or not)
{\em in cohomological degree one}. As a prototype: let $M$ be a compact K\"ahler
manifold with fundamental group $\Ga$. For all $i\in\N$, let $b^i(\G)=\rk H^i(\G,\Z)={\rm dim}_\R\, H^i(\G,\R)={\rm dim}_\C\, H^i(\G,\C)$. Classical Hodge theory shows that $b^1(\Ga)=b^1(M)$ has to be even. Considering finite
\'etale covers one obtains that $b^1(\Ga')$ has to be even for any finite index
subgroup $\Ga'$ of $\Ga$. For example, the free group on $n$-generator
$F_n$ is never K\"ahler, $n\geq 1$.

\begin{rema}
In \cite{Voisin}, Voisin proved that there exist compact K\"ahler
manifolds whose homotopy type cannot be realized by a smooth
projective manifold. However, the question whether or not any K\"ahler group
is a projective group  is still open.
\end{rema}

A very interesting conjecture concerning infinite
K\"ahler groups, due to Carlson and Toledo and publicized by Koll\'ar \cite{ko}, deals with {\em cohomology in degree $2$}:
\begin{conj}[Carlson-Toledo] \label{conj1}
Let $\Ga$ be an infinite K{\"a}hler group. Then virtually $H^2(\G,\R)\neq 0$.
\end{conj}

We refer to Section~\ref{explanation} for different interpretations of
this conjecture. Notice that it holds trivially for fundamental groups
of complex curves.
Although K\"ahler groups are not necessarily residually finite (and
thus not necessarily linear) \cite{to}, all known infinite such groups admit a linear representation with
{\em unbounded image}. We will from now on restrict to this case.

\sspace
The strongest evidence for Carlson-Toledo's conjecture is the
following theorem. The finite dimensional case is folkloric (we give a proof in Section~\ref{sec2} for the sake of
completeness) and the Hilbert space case is one of the main results of Reznikov in~\cite{Re}:
\begin{theo} \label{prop2}
Let $\Ga$ be a K\"ahler group. If $\Ga$ admits a linear representation $\rho:
\Ga \lo G$, with $G$ the linear group ${\rm GL}(V)$ of a finite dimensional
complex vector space $V$ or the isometry group of a Hilbert space, satisfying
$H^1(\Ga, \rho) \not = 0$,
then $b^2(\Ga) >0$.
\end{theo}

Recall that a finitely generated group is said {\em schematically rigid} if all linear representations of $\G$ are {\em rigid}, namely if for any $n \in \nn$
and any representation $\rho: \Ga \lo \GL(n, \cl)$ one has $H^1(\Ga,
\Ad \rho)=0$. As a corollary to Theorem~\ref{prop2} one has:

\begin{coro}
Let $\Ga$ be an infinite K\"ahler group. If $\Ga$ does not satisfy
Carlson-Toledo's conjecture then necessarily:
\begin{itemize}
\item[(a)] $\Ga$ has Kazhdan's property $(T)$.
\item[(b)] $\Ga$ is schematically rigid. 
\end{itemize}
\end{coro}

Let $M$ be a compact K\"ahler manifold and assume that $\Ga=\pi_1(M)$ is schematically rigid
and admits a reductive linear representation $\rho: \Ga \lo
\GL(n, \C)$. By Simpson's theory \cite{Sim}, the representation $\rho$ is then the monodromy of a (polarized) complex variation of Hodge structure ($\C$-VHS for short) over $M$, see Section~\ref{var}.  
The first result of this paper states that
Conjecture~\ref{conj1} holds under some
assumption on the variation of Hodge structure:
\begin{theo}\label{cartoun}
Let $M$ be a compact K\"ahler manifold. Assume that there exists an
irreducible rigid unbounded reductive representation $\rho:\pi_1(M)\fd
\GL(n,\C)$ for some $n>0$ and that the variation of Hodge
structure $(E,\t)=(\oplus_{i=0}^k  E^i,\oplus_{i=0}^{k-1}\t_i)$ over $M$ which induces $\rho$ satisfies $\rk E^i=1$ for some $0<i<k$. Then $H^2(\pi_1(M),\R)\not=0$.
\end{theo}

\begin{rema}
The theorem holds under the more general assumption that the rank of $\t_i:T_M\fd\Hom(E^i,E^{i+1})$ is $1$ for some $0\leq i\leq k-1$ (see the proof in Section~\ref{proofofcartoun}). 
\end{rema}

Note in particular that this result (and the discussion in section~\ref{conj}) implies
that the Carlson-Toledo conjecture holds whenever the image
of the rigid representation sits in some ${\rm U}(p,1)$.\\

\sspace

The fact that the representation $\rho$ comes from a $\C$-VHS implies in particular that $\rho(\G)$ sits in some $\U(p,q)$, $p+q=n$. Moreover there is a period domain $D$, which is an open orbit of the group $\PU(p,q)$ in a manifold of flags of $\C^n$, and a period map, which is a $\rho$-equivariant holomorphic superhorizontal map from the universal cover of $M$ to $D$. We call such period domain and map the $\PU(p,q)$ period domain and map associated to the $\C$-VHS, see Section~\ref{var} for details. It turns out that the Carlson-Toledo conjecture can be related to the question of whether the period map kills the second homotopy group $\pi_2(M)$ of $M$ or not. Indeed if the answer is yes, then the conjecture is true. Since the period map is superhorizontal, a first thing to study is the existence of non homotopically trivial superhorizontal 2-spheres in the period domain $D$. In Section~\ref{superhori}, using Gromov's $h$-principle
we prove the following result, which is of independent geometric interest:
\begin{theo}\label{horiz}
Let $(E,\t)$ be a $\C$-VHS over a K\"ahler manifold and let
$D$ be the $\PU(p,q)$ period domain associated to $E$. Suppose $(E,\t)$ does not
satisfy the hypotheses of Theorem~\ref{cartoun}. Then $\pi_2(D)$ can be generated by superhorizontal 2-spheres. 
\end{theo}

It is also customary to consider ``smaller'' period domains that one can associate to a $\C$-VHS: the real Zariski-closure $G_0$ of $\rho(\Ga)$ in $\GL(n,\C)$ is a group of Hodge type by~\cite{Sim} and the period map can be chosen to land in a period domain which is an open $G_0$-orbit in some flag manifold $G/Q$, where $G$ is the complex Zariski-closure
of $\rho(\G)$ and $Q$ a parabolic subgroup of $G$, see Sections~\ref{perdo} and~\ref{var}. We will comment on the generalization of Theorem~\ref{horiz} to such period domains at the end of the paper.   

\medskip

{\em Acknowledgements.} We wish to thank the referee for his helpful comments and for having suggested many improvements in the exposition.  

\section{What does Conjecture~\ref{conj1} mean?} \label{explanation}

\subsection{In terms of group extensions}
Recall that for a group $\G$ and an Abelian group $A$ with trivial
$\Gamma$-module structure, the group $H^2(\Ga, A)$ classifies the
central $A$-extensions
$$ 0 \lo A \lo \tilde{\Gamma} \lo \Gamma \lo 1$$
of $\Ga$. As $b^2(\Ga) = \textnormal{rk}\, H^2(\Ga, \z)$,
Conjecture~\ref{conj1} means that any infinite K\"ahler group admits
(after maybe passing to a finite index subgroup) a non-trivial
central $\z$-extension (which does not trivialize when restricted to
any finite index subgroup).

\subsection{In terms of classical topology}

For any group $\Ga$ the universal coefficients exact
sequence yields the isomorphism
$$H^2(\Ga, \real) = \Hom_\real(H_2(\Ga, \real), \real).$$
In particular it is equivalent to show $b_2(\Ga) >0$ or $b^2(\Ga)>0$.

\sspace
For any reasonable topological space $M$ with fundamental group $\Ga$
the universal cover $p:\tilde{M}\fd M$ is a principal $\Ga$-bundle over
$M$. Thus it defines (uniquely in the homotopy category) a morphism $c: M \lo B\Ga$
from $M$ to the classifying Eilenberg-MacLane space $B\Ga =
K(\Ga,1)$. The induced morphism 
$c_{\star}:H_\star(M, \R) \lo H_\star(\Ga, \R)$ is easily seen to be an
isomorphism in degree $1$ and an epimorphism in degree $2$: 
$$
H_2(M, \R) \twoheadlongrightarrow H_2(\Ga, \R).
$$
Dually:
$$
H^2(\Ga, \R) \hooklongrightarrow H^2(M, \R).
$$

\sspace
How can we characterize the quotient $H_2(\Ga, \R)$ of $H_2(M, \R)$?
In fact this quotient first appeared in Hopf's work on the Hurewicz
morphism comparing homotopy and homology:
\begin{theo} [Hopf] \label{Hopf}
Let $N$ be a paracompact topological space. Let $c: N \lo B\pi_1(N)$
be the classifying morphism and $h: \pi_\star(N) \lo H_\star(N, \z)$ the
classical Hurewicz morphism.
Then the sequence of Abelian groups
$$
\pi_2(N) \stackrel{h}{\lo} H_2(N, \z) \stackrel{c_\star}{\lo} H_2(\pi_1(N),
\z) \lo 0
$$
is exact.
\end{theo} 

Cohomologically:
\begin{coro} \label{corol3}
Let $N$ be a paracompact topological space and $\pi_2(N)
\otimes_\z \R \stackrel{h}{\lo} H_2(N, \R)$ the Hurewicz
morphism. Then:
$$H^2(\pi_1(N), \R) = \{[\om] \in H^2(N, \R)\; |\; \forall \phi:S^2
\lo N, \; <[\om], \phi_\star [S^2]>=0\} \subset H^2(N, \R)\;\;,$$ where 
$<\cdot, \cdot>: H^2(N, \R) \times H_2(N, \R) \lo \R$ is the
natural non-degenerate pairing between homology and cohomology.
\end{coro}

\begin{rema}
Nowadays Theorem~\ref{Hopf} is a direct
application of the Leray-Cartan spectral sequence.
\end{rema}

Since $\pi_2(M)$ is nothing else than $H_2(\tilde{M}, \z)$, Carlson-Toledo's conjecture can be restated:
\begin{conj} \label{conj2}
Let $M$ be a compact K\"ahler manifold with infinite fundamental group. There exists a finite \'etale cover $M'$ of $M$ such that the natural morphism $H_2(\tilde{M}, \real) \fd H_2(M', \real)$ is not surjective, or, equivalently, such that the natural morphism
$H^2(M', \real) \fd H^2(\tilde{M}, \real)$ is not
injective. 
\end{conj}


Note that for any compact manifold $M$ and any finite \'etale cover
$M'$ of $M$ the arrow $H^2(M, \real) \lo H^2(M', \real)$ is injective
by the projection formula.

\subsection{In terms of $\cl^\star$-bundles}
Recall that for a reasonable topological space $M$ the group $H^2(M,
\z)$ canonically identifies with the group 
$\mathcal{L}(M)$ of principal $\cl^\star$-bundles: on the one hand
$H^2(M, \z) = [M, K(\z, 2)]$, on the other hand 
$\mathcal{L}(M)= [M, B\cl^\star]$. But both $K(\z, 2)$ and $B\cl^\star$ have
as canonical model the infinite projective space $\cl
\proj^{\infty}$. Thus Conjecture~\ref{conj1} states that any
infinite K\"ahler group $\Ga$ admits a finite index subgroup $\Ga'$
whose classifying space $B\Ga'$ supports a non-trivial
principal $\cl^\star$-bundle. Let $M$ be a compact K\"ahler manifold 
with infinite fundamental group $\Ga$. As $\tilde{M}$ is the homotopy
fiber of $M \lo B\Gamma$, Conjecture~\ref{conj2}
says there exists a finite \'etale cover $M'$ of $M$ and a non-trivial
principal $\cl^\star$-bundle on $M'$ whose pull-back to $\tilde{M}$ becomes
trivial. In these statements ``non-trivial'' means ``with non-trivial
{\em rational} first Chern class''.

\subsection{Equivalence} \label{equivalence}
The equivalence of these three points of view is clear. Given $M$ and a line bundle $L$ on $M$ with associated
principal $\cl^\star$-bundle $L^0$, one can consider the long homotopy exact
sequence for the fibration $L^0 \lo M$:
\begin{equation} \label{long}
 \cdots \lo \pi_2(L^0) \lo \pi_2(M) \stackrel{\partial_L}{\lo}\pi_1(\cl^\star) = \z \lo
\pi_1(L^0) \lo \pi_1(M) \lo 1 \;\;.
\end{equation}
The boundary map $\partial_L: \pi_2(M) \lo \z$ is just the first Chern class
map of $L$ restricted to $\pi_2(M)$: if $[\alpha] \in \pi_2(M)$ is
represented by $\alpha: S^2 \lo M$ then $\partial_L([\alpha])=
<\alpha^\star (c_1(L)), [S^2]> \in \z$. As $H_2(\tilde M, \z) = \pi_2(M)$ and
$H^2(\tilde{M}, \real)$ is dual to $H_2(\tilde M, \real)$ the principal $\cl^\star$-bundle
$p^\star L^0$ is trivial if and only if $\partial_L:
\pi_2(M) \lo \real$ is zero. Then the long exact sequence~(\ref{long})
gives the short exact sequence: 
$$ 0 \lo \z \lo \pi_1(L^0) \lo \pi_1(M) \lo 1 \;\;,$$
the element of $H^2(\pi_1(M),\Z)$ classifying this extension being the first Chern class of $L$.

\section{Carlson-Toledo's conjecture: the non-rigid case} \label{sec2}

\subsection{The reductive case}
The simplest instance of Theorem~\ref{prop2} is the following:
\begin{lemma} \label{lef}
Let $\Ga$ be a K\"ahler group. If $b^1(\Ga)>0$ then $b^2(\Ga)>0$.
\end{lemma}

Essentially the same proof generalizes to:

\begin{lemma}\label{l1}
{\rm (\cite{Re})}. Let $\Ga$ be a K\"ahler group. If $\Ga$ admits a linear representation $\rho:
\Ga \lo G$ with Zariski-dense image, where $G$ denotes a linear
{\em reductive} complex algebraic group, and satisfying $H^1(\Ga, \rho) \not = 0$, 
then $b^2(\Ga) >0$.
\end{lemma}

\begin{demo}
Let $M$ be a compact K\"ahler manifold whose fundamental group is isomorphic to $\Ga$. Let $E_{\rho}
\lo M$ be the complex local system on $M$ associated with $\rho$.
Consider the following commutative diagram:
$$
\xymatrix{
H^{1}(\Ga, \rho) \times H^{1}(\Ga, \rho^{\star}) \ar[r] \ar@{=}[d]  &H^{2}(\Ga,
\cl) \ar[r] \ar@{^{(}->}[d] &\cl \ar@{=}[d]\\
H^{1}_{DR}(M, E_{\rho}) \times H^{1}_{DR}(M, E_{\rho}^{\star}) \ar[r] &H^{2}_{DR}(M, \cl)
\ar[r] &\cl \\
}
$$
where $\rho^\star$ denotes the contragredient representation of $\rho$ and the second row
of the diagram is given by
$$(\alpha, \beta) \lo \alpha \wedge \beta \lo \int_{M} \alpha \wedge \beta\wedge \omega_{M}^{{\rm dim}\,M-1}\;\;.$$ 
By the hard Lefschetz theorem generalized by Simpson to the case of
reductive representations \cite[Lemmas 2.5, 2.6]{Sim} this second row
defines a non-degenerate bilinear form. This implies the result.
\end{demo}

\subsection{The general linear case} \label{l2}

\begin{lemma} \label{reduc}
Let $\Ga$ be a K\"ahler group. If $\Ga$ admits a linear representation
$\rho: \Ga \lo {\rm GL}(n, \cl)$ satisfying $H^1(\Ga, \rho) \not = 0$ 
then $b^2(\Ga) >0$.
\end{lemma}

\begin{demo}
Let $G= \RR \ltimes \UU$ be the Levi decomposition of the
complex Zariski-closure $G$ of $\rho(\Ga)$ in
${\rm{GL}}(n,\cl)$, where $\UU$ denotes the unipotent radical of
$G$ and $\RR$ a Levi factor lifting the reductive quotient $G
/\UU$. Let 
$l(\rho) : \Ga \lo \RR$ be the quotient representation of
$\rho$. As $\rho$ has Zariski-dense image in $G$, the representation
$l(\rho)$ has Zariski-dense image in the reductive group $\RR$.

\sspace
If $\UU$ is non-trivial let $$\UU = \UU_1
\supset \cdots \supset \UU_i \supset \cdots \supset  \UU_r$$  
be the central descending series for $\UU$ (i.e.
$\UU_1=\UU$, $\UU_{i+1}= \overline{[\UU, \UU_{i}]}$ is the
Zariski-closure of the subgroup of $\UU$ generated by the commutators
of $\UU$ and $\UU_{i}$). The $\UU_i$'s are unipotent, the quotients 
$\A_i:=\UU_i/\UU_{i+1}$ are non-trivial of additive type. The action of $\RR$
on $\UU$ by automorphisms preserves the central descending
series. Let us define $G_i=\RR \ltimes \A_i$ and $\tau_i \in
H^{1}(G_i, \A_i)$ be the canonical class of the semi-direct
product $G_i$.

\sspace
Let $\rho_1 : \Ga \lo \GG_1$ be the quotient representation of $\rho$.
If the class $\tau_{1, \Ga}= \rho_{1}^{\star}(\tau_{1}) \in H^{1}(\Ga,
(\A_1)_{l(\rho)})$ does not vanish, the conclusion  $H^{2}(\Ga,
\cl)\not = 0$ follows from Lemma~\ref{l1} applied to the reductive
representation  $(\A_1)_{l(\rho)}$.
If $\tau_{1, \Ga}=0$ then up to conjugacy by an element in $\GG$
one can assume that the image $\rho_{1}(\Ga)$ is contained in
$\RR$. Thus the image $\rho(\Ga)$ is contained in $\RR\ltimes \UU_2$.
The quotient representation $\rho_2: \Ga \lo \GG_2$ of $\rho$ is
then well-defined, thus also the class $\tau_{2,\Ga}= \rho_{2}^{\star}(\tau_{2}) \in H^{1}(\Ga,
(\A_2)_{l(\rho)})$.

\sspace
By induction, for the first step $i$ where the class
$\tau_{i,\Ga}$ does not vanish one deduces $H^{2}(\Ga,
\cl)\not = 0$ from Lemma~\ref{l1} applied to the reductive
representation $(\A_i)_{l(\rho)}$. Thus one can assume that all the
classes $\tau_{i,\Ga}$, $1 \leq i \leq r$, vanish. Hence $\rho(\Ga) \subset \RR$, $\UU=\{1\}$ and
$\GG=\RR$. But then $\rho$ is reductive and satisfies by hypothesis $H^1(\Ga, \rho) \not = 0$, thus $H^{2}(\Ga,\cl)\not = 0$ by Lemma~\ref{l1} applied to $\rho$.
\end{demo}

\subsection{The Hilbert case}
The proof of Lemma~\ref{lef} generalizes also to the case where $\rho$
is a unitary representation in some Hilbert space. This is the main
result of Reznikov in
\cite{Re}, to which we refer for a proof.

\section{Carlson-Toledo conjecture and variations of Hodge structure}\label{back}

\subsection{Period domains}\label{perdo}

We recall the necessary background about period domains. We refer to \cite{GS} for a detailed account of the subject.

Let $G$ be a connected complex semi-simple Lie group and $G_0$ a noncompact real form of $G$. We assume that $G_0$ is a group of Hodge type, namely that it has a compact maximal Abelian subgroup $T_0$ (see~\cite[4.4.4]{Sim}). Let $Q$ be a parabolic subgroup of $G$ which contains $T_0$ and such that $V_0:=G_0\cap Q$ is compact. 

We choose a maximal compact subgroup $K_0$ of $G_0$ which contains $V_0$, and a maximal compact subgroup $G_u$ of $G$ which contains $K_0$.
The Lie algebras of $G$, $G_0$ and $Q$ will be denoted by $\gg$, $\gg_0$ and $\qg$. Let $B$ be the Killing form on $\gg$. The Lie algebra $\gg_u$ of $G_u$ is a real form of $\gg$ and we denote complex conjugation with respect to $\gg_u$ by $\tau$. Let $\ng$ be the maximal nilpotent ideal of $\qg$ and $\vg=\qg\cap\tau(\qg)$. Then $\qg=\vg\oplus\ng$ and $\gg=\vg\oplus\ng\oplus\tau(\ng)$. The Lie algebra of $V_0$ is $\vg_0=\gg_0\cap\qg=\gg_0\cap\qg\cap\tau(\qg)=\gg_0\cap\vg$. We shall denote the Lie algebra of $K_0$ by $\kg_0$ and its complexification by  $\kg$. We have $\vg\subset \kg$.

The complex analytic quotient space $D_u=G/Q$ is called a {\em flag manifold}. The group $G$ is a holomorphic principal $Q$-bundle over $D_u$ and the holomorphic tangent bundle of $D_u$ is the bundle associated to this principal bundle by the adjoint action of $Q$ on $\gg/\qg\simeq\tau(\ng)$. The $G_u$-orbit of $eQ\in D_u=G/Q$ is open because the real span of $\gg_u$ and $\qg$ is all of $\gg$ and it is closed, since $G_u$ is compact. Hence $G_u$ acts transitively on $D_u$ with isotropy group $V_u=G_u\cap Q=V_0$ and we identify the quotient space $G_u/V_u$ with $D_u$. 

Open orbits of noncompact real forms of $G$ in the flag manifold $G/Q$ are called {\em flag domains} (\cite{FHW}). With our assumptions, this is the case of the $G_0$-orbit $D$ of $eQ$, which can be identified with the quotient space $G_0/V_0$.  Because $G_0$ is of Hodge type and $V_0=G_0\cap Q$ is compact, we will call the homogeneous complex manifold $D$ a {\em period domain}. We say that $D$ is ``dual'' to $D_u$.

Note that any parabolic subgroup $Q'$ of $G$ such that $V_0=G_0\cap Q'$ induces a complex structure on $D$.

As a $C^\infty$ vector bundle, the complex tangent bundle $T_{D_u}$ (resp. $T_{D}$) to $D_u$ (resp. $D$) is associated to the principal $V_0$-bundle $G_u$ (resp. $G_0$) by the adjoint representation of $V_0$ on $\gg/\qg$.
We endow the tangent space of $D_u$ at $eV_0$, which is naturally isomorphic to $\tau(\ng)$, with the ${\rm Ad} V_0$-invariant inner product $(\xi,\zeta)\mapsto-B(\xi,\tau(\zeta))$. By $G_u$-translation, this inner product gives rise to a $G_u$-invariant Hermitian metric $g_{D_u}$ on $D_u$ and by $G_0$-translation, to a $G_0$-invariant Hermitian metric $g_{D}$ on $D$.

Let $X=G_0/K_0$ be the symmetric space associated to $G_0$ and $X_u=G_u/K_0$ its compact dual. Then, we have natural fibrations $\pi_u:D_u=G_u/V_u\fd G_u/K_0=X_u$ and $\pi:D=G_0/V_0\fd G_0/K_0=X$. The fibers of $\pi_u$, resp. $\pi$, are complex submanifolds of $D_u$, resp. $D$, all isomorphic to the central fiber $K/K\cap Q=K_0/V_0$ (in both cases, although $\pi\not={\pi_u}_{|D}$) and the fibration is $G_u$-invariant, resp. $G_0$-invariant. The vectors in $T_{D}$ which are tangent to the fibers of $\pi$ form a $C^\infty$ subbundle $T^v_{D}$. A vector $\xi$ in the fiber of $T_{D}$ over $eV$, which is identified with $\gg/\qg$, is tangent to the fiber of $\pi$ if and only if $\xi\in\kg/\kg\cap\qg$. It follows that $T_{D}^v$ is associated to the $V_0$-bundle $G_0$ by the adjoint representation of $V_0$ on $\kg/\kg\cap\qg$.

The maximal compact subgroup $K_0$ of $G_0$ determines a Cartan decomposition $\gg=\kg\oplus\pg$. The adjoint action of $V_0$ on $\pg/\pg\cap\qg$ associates a $C^\infty$ vector bundle $T^h_{D}$ to the $V_0$-principal bundle $G_0$ and since $\gg/\qg=\kg/\kg\cap\qg\oplus\pg/\pg\cap\qg$, $T_{D}=T^v_{D}\oplus T^h_{D}$ is a $G_0$-invariant orthogonal splitting of the tangent bundle into two $C^\infty$ subbundles. Moreover, if we endow $X$ with its natural $G_0$-invariant metric then $\pi$ is a Riemannian submersion, meaning that $d\pi_{|T^h_{D}}:T^h_{D}\fd T_{X}$ is an isometry.

Of course, we can make the same kind of constructions on $D_u$ (with obvious changes) and we get in particular a $G_u$-invariant decomposition $T_{D_u}=T^v_{D_u}\oplus T^h_{D_u}$.

The Hermitian metric $g_{D_u}$ (resp. $g_{D}$) induces a $G_u$-invariant metric $h_{D_u}$ (resp. a $G_0$-invariant metric $h_{D}$) on the canonical bundle $K_{D_u}=\Lambda^{\dim_\C D_u} T^\star_{D_u}$ of $D_u$ (resp. the canonical bundle $K_{D}={K_{D_u}}_{|D}$ of $D$). The computations of \cite{GS} show that

\begin{prop}\label{cano}
The $G_u$-invariant curvature form $\Th_{D_u}(K_{D_u})$ of $(K_{D_u},h_{D_u})$ is negative, and the $G_0$-invariant curvature form $\Th_{D}(K_{D})$ of $(K_{D},h_{D})$ is positive on the horizontal space $T^h_{D}$ and negative on the vertical space $T^v_{D}$. In particular, $\log(h_{D_u}/h_{D})$ is a well-defined function on $D$ whose Levi form is positive definite on $T^h_{D}$.
\end{prop}

\sspace

It follows from~\cite[Chap. 4.A]{BR} that there exists an element $\xi$ in the center of $\vg$ such that all the eigenvalues of $\ad\xi$ acting on $\gg$ belong to $\sqrt{-1}\,\Z$ and moreover, if $\gg_\ell$ denotes the eigenspace associated to the eigenvalue $\sqrt{-1}\,\ell$, then if $\ng^{(r)}$ denotes the $r$-th step of the central descending series of $\ng$,  
$$
\qg=\bigoplus_{\ell\geq 0}\gg_\ell\ \mbox{ and }\ \ng^{(r)}=\bigoplus_{\ell\geq r}\gg_\ell.
$$
For all $m>0$, let $W_m=\oplus_{\ell\geq -m}\,\gg_\ell$. Then $W_m/\qg\subset \gg/\qg$ is an ${\rm Ad}\,Q$-invariant subspace, hence defines a holomorphic subbundle of $T_{D_u}$ (or $T_{D}$). In particular, for $m=1$, we get a holomorphic subbundle $H_{D_u}$, resp. $H_{D}={H_{D_u}}_{|D}$, of $T_{D_u}$, resp. $T_{D}$, which is included in $T^h_{D_u}$, resp. $T^h_{D}$, and called the {\em superhorizontal distribution} of $D_u$, resp. $D$. Moreover, $\gg_{1}$ Lie-generates $\ng$ and this implies that the superhorizontal distribution satisfies the so-called bracket-generating condition (see for example~\cite{Pansu}).

\subsection{Measurable open orbits in flag manifolds}\label{mes}

We keep the assumptions and notation of the previous section, and we consider a new parabolic subgroup $Q'$ of $G$ which also contains the maximal Abelian subgroup $T_0$ of $G_0$, but this time we do not assume that $G_0\cap Q'$ is compact. To emphasize this difference, we set $L_0=G_0\cap Q'$.

It follows from \cite[Sect. 4.2 and 4.5]{FHW} that the $G_0$-orbit $D'$ of $eQ'$ in $D'_u$ is open and measurable, meaning that it admits a $G_0$-invariant, possibly indefinite, K\"ahler metric. We will identify $D'$ with $G_0/L_0$ and $D'_u$ with $G_u/L_u$ where $L_u$ is the compact real Lie subgroup $G_u\cap Q'\subset G$. $D'$ contains a special compact complex submanifold $C\simeq K_0/K_0\cap L_0$ (\cite[Sect. 4.3]{FHW}).

There is a decomposition $\qg'=\lg\oplus\ng'$, where $\lg$ is the complexification of the Lie algebra of $L_0$ and $\ng'$ is the maximal nilpotent ideal of $\qg'$. As in the case of period domains, $G$ can be seen as a holomorphic principal $Q'$-bundle on $D'_u$ and the holomorphic tangent bundle of $D'_u$ is the bundle associated to this principal bundle by the adjoint action of $Q'$ on $\gg/\qg' \simeq \tau(\ng')$. Let us mention that the underlying $C^\infty$ vector bundle on $D'$, which is associated to the principal $L_0$-bundle $G_0$, admits in general no decomposition similar to the decomposition in horizontal and vertical components. Nevertheless, we can endow the canonical bundle $K_{D'_u}$ of $D'_u$ (resp. $K_{D'}$ of $D'$) with a $G_u$-invariant (resp. $G_0$-invariant) metric $h_{D'_u}$ (resp. $h_{D'}$) and an analogous statement to Proposition~\ref{cano} (see~\cite[Sect. 4.7]{FHW}) is available:

\begin{prop}\label{canon}
The $G_u$-invariant curvature form $\Th_{D'_u}(K_{D'_u})$ of $(K_{D'_u},h_{D'_u})$ is negative, and the $G_0$-invariant curvature form $\Th_{D'}(K_{D'})$ of $(K_{D'},h_{D'})$ has signature $(\dim_\C D'-\dim_\C C,\dim_\C C)$. In particular, $\varphi:=\log(h_{D'_u}/h_{D'})$ is a well-defined function on $D'$ whose Levi form ${\cal L}(\varphi)$ has everywhere at least $\dim_\C D'-\dim_\C C$ positive eigenvalues.
\end{prop}

We will need later a more precise description of the tangent vectors on which ${\cal L}(\varphi)$ is positive definite. Recall that the tangent bundle of $D'$ is the bundle associated to the principal $L_0$-bundle $G_0$ by the adjoint representation of $L_0$ on $\gg/\qg'$. From the computations of \cite{FHW}, $\Th_{D'}(K_{D'})$ is positive on any tangent vector $(e,\xi)\in G_0\times_{L_0} \gg/\qg'\simeq T_{D'}$ such that $\xi\in\pg/\pg\cap\qg'$. By $G_0$-invariance of $\Th_{D'}(K_{D'})$, it follows that it is positive on any tangent vector to $D'$ which has a representative $(g,\xi)$ with $\xi\in\pg/\pg\cap\qg'$ and then, on such a tangent vector, ${\cal L}(\varphi)$ is also positive.

\subsection{Complex variations of Hodge structure and $\PU(p,q)$ period domains}\label{var}

We refer to \cite{Gri} and \cite{Sim} for details on this section.

\begin{defi}
Let $M$ be a complex analytic manifold. A (polarized) complex variation of Hodge structure ($\C$-VHS for short) on $M$ is a $C^\infty$ complex vector bundle $E$ (of finite rank) with a decomposition $E=\oplus_{i=0}^k E^i$ in a direct sum of $C^\infty$ subbundles, a flat connection $\nabla$ satisfying Griffiths's transversality condition
$$\nabla:E^i\fd A^{0,1}(E^{i-1})\oplus A^{1,0}(E^i)\oplus A^{0,1}(E^i)\oplus A^{1,0}(E^{i+1})
$$
and a parallel Hermitian form $h$ which makes the Hodge decomposition orthogonal and which on $E^i$ is positive if $i$ is even and negative if $i$ is odd. If we set $r_i=\rk E^i$, the monodromy representation $\rho$ of the flat connection $\nabla$ has image contained in a group isomorphic to $\U(p,q)\subset\GL(p+q,\C)$ where $p=\sum_{i\ {\rm even}}r_i$ and $q=\sum_{i\ {\rm odd}}r_i$.
\end{defi}

Assume we have a complex variation of Hodge structure on a complex manifold $M$. We have the so-called Hodge filtration of $E$ by the subbundles $F^s=\oplus_{i\leq s} E^i$ which are holomorphic for the connection $\nabla$.

Let us fix $x_0\in M$, and let $\E=\oplus_{i=0}^k\E^i$ be the fiber of $E$ over $x_0$ endowed with the restriction of the Hermitian metric $h$ on $E$, which will be denoted by the same letter. Then, using notation which are consistent with those in Section~\ref{perdo}, we may set $G=\PGL(\E)$, $G_0=\PU(\E,h)$, $G_u=\PU(\E)$ (w.r.t. the standard Hermitian form on $\E$), and let $Q$ be the parabolic subgroup of $G$ which stabilizes the flag 
${\mathbb F}^0\subset {\mathbb F}^{1}\subset\dots\subset {\mathbb F}^{k-1}\subset{\mathbb F}^k=\E$ defined by the Hodge filtration and $V_0=V_u=G_0\cap Q={\rm P}(\U(\E^0,{h}_{|\E^0})\times\dots\times\U(\E^k,{h}_{|\E^k}))$. We obtain a period domain $D=G_0/V_0$ in the flag manifold $D_u=G/Q=G_u/V_u$. We call $D$ the $\PU(p,q)$ period domain associated to the $\C$-VHS $E$. 

Choose bases of each of the $\E^i$ and put them together in order to get a basis of $\E$. Identifying now $\gg=\sg\lg(\E)$ with $\sg\lg(p+q,\C)$, the usual Cartan subalgebra $\hg\subset\gg$ is the set of (tracefree) diagonal matrices. We denote by $\Delta$ the root system associated to $\hg$  and $\Delta^+$ the set of positive roots with respect to the usual ordering of $\Delta$. Then, $\Pi=\{\a_1,\dots,\a_{p+q-1}\}\subset\Delta$ (where $\a_j({\rm diag}(h_1,\dots,h_{p+q}))=h_{j+1}-h_j)$ is a simple system. For each root $\a\in\hg^\star$, we denote by $\gg^\a$ the corresponding root-space. The parabolic subalgebra $\qg$ is associated to the subset $\Pi_\qg=\{\a_{r_0},\a_{r_0+r_1},\dots,\a_{r_0+\dots+r_{k-1}}\}\subset\Pi$ in the sense that if $\Phi=\Delta^+\cup\{\a\in\Delta\,|\,\a\in{\rm span}(\Pi\backslash\Pi_\qg)\}$ then
$$\qg=\hg\oplus\bigoplus_{\a\in\Phi}\gg^\a.
$$
We also have
$$\vg=\hg\oplus\bigoplus_{\a\in\Phi\cap-\Phi}\gg^\a\ \ {\rm and}\ \ \ng=\bigoplus_{\substack{\a\in\Phi\\ \a\not\in-\Phi}} \gg^\a.
$$
Moreover, in this case, $\gg_\ell=\oplus_{i-j=\ell}\Hom(\E^i,\E^j)\cap \sg\lg(\E)$ for all $\ell$, so that $\qg=\oplus_{i=0}^k\Hom(\E^i,{\mathbb F}^{i})\cap \sg\lg(\E)$ and 
the superhorizontal distribution is the holomorphic bundle associated to $(\oplus_{i=0}^{k-1}\Hom(\E^i,{\mathbb F}^{i+1})\cap \sg\lg(\E))/\qg$ by the adjoint action of $Q$.

Let $K_0={\rm P}(\U(\oplus_{i} \E^{2i},{h}_{|\oplus_{i} \E^{2i}})\times \U(\oplus_{i} \E^{2i+1},{h}_{|\oplus_{i} \E^{2i+1}}))$. This is a maximal compact subgroup of $G_0$ containing $V_0$. The symmetric space $X=G_0/K_0$ is Hermitian but we stress that $\pi:D\fd X$ is never holomorphic if $D$ is endowed with the complex structure defined by $Q$ (except if $V_0=K_0$).

\medskip
Geometrically, ${\rm PU}(p,q)$-period domains and flag manifolds can be described as follows. We endow ${\mathbb E}=\C^{p+q}$ with the Hermitian form $h=\sum_{i=1}^p|z_i|^2-\sum_{i=1}^q|z_{p+i}|^2$. The flag manifold $D_u$ is simply the space of flags $\{0={F}^{-1}\subset{F}^0\subset {F}^{1}\subset\dots\subset {F}^{k-1}\subset {F}^k=\C^{p+q}\}$ of $\C^{p+q}$ with $\dim F^{i}/F^{i-1}=r_i$, $p=\sum_{i\ {\rm even}}r_i$ and $q=\sum_{i\ {\rm odd}}r_i$. 
The period domain $D$ is the open set of $D_u$ consisting of those flags such that for all $-1\leq i\leq k-1$, $(-1)^{i} h$ is negative definite on ${F^i}^\perp$, the orthogonal of $F^i$ in $F^{i+1}$ with respect to the indefinite form $h$. The symmetric space $X$ is the space of $p$-planes in $\C^{p+q}$ on which $h$ is positive definite. It is an open set in $X_u$, the Grassmannian of $p$-planes in $\C^{p+q}$. The fibration $\pi:D\fd X$ is given by 
$$
\{0={F}^{-1}\subset{F}^0\subset {F}^{1}\subset\dots\subset {F}^{k-1}\subset {F}^k=\C^{p+q}\}\longmapsto\bigoplus_{i\ {\rm odd}}{F^i}^\perp,
$$ 
whereas $\pi_u:D_u\fd X_u$ is given by 
$$
\{0={F}^{-1}\subset{F}^0\subset {F}^{1}\subset\dots\subset {F}^{k-1}\subset {F}^k=\C^{p+q}\}\longmapsto\bigoplus_{i\ {\rm odd}}{F^i}^{\perp_u},
$$ 
where ${F^i}^{\perp_u}$ is the orthogonal of $F^i$ in $F^{i+1}$, with respect to the standard definite Hermitian form $\sum_{i=1}^{p+q}|z_i|^2$ on $\C^{p+q}$. Remark that the central fibers of $\pi$ and $\pi_u$, namely the fibers above the $p$-plane 
$\{z_{p+1}=\dots =z_{p+q}=0\}$, are the same.

\medskip

Now, let $\tilde M$ be the universal cover of $M$. The complex variation of Hodge structure $E$ induces a $\rho$-equivariant map $f:\tilde M\fd D$ called the {\em period map}, and the $\C$-VHS $E$ is said to be nontrivial if $f$ is a nonconstant map. The transversality condition implies that $f$ is holomorphic with respect to the complex structure induced by $Q$ and moreover that $df(T_{\tilde M})\subset H_{D}$: we say that $f$ is a holomorphic {\em superhorizontal map}. 

The link between the $\C$-VHS and the period map $f$ can be seen as follows. The pull-back 
$f^\star H_{D}$ is a holomorphic bundle on $\tilde M$ which goes down on $M$ by $\rho$-equivariance. As such, it can be identified with $\oplus_{i=0}^{k-1}{\rm Hom}(E^i,E^{i+1})$ where the vector bundles $E^i\simeq F^i/F^{i-1}$ are endowed with the holomorphic quotient structure. The $(1,0)$ part of the differential of $f$ can hence be interpreted as a morphism $T_M\fd\oplus_{i=0}^{k-1}{\rm Hom}(E^i,E^{i+1})$. On the other hand, the flat connection $\nabla$ on $E$ gives rise to ${\cal O}_M$-linear maps $\theta_i:E^i\fd E^{i+1}\otimes\Omega^1_M$ which are holomorphic sections of ${\rm Hom}(E^i,E^{i+1})\otimes\Omega^1_M$. The so-called Higgs field $\theta:=\oplus_{i=0}^{k-1} \theta_i$ is therefore another incarnation of $df$. Note that $\t$ satisfies the relation $[\t,\t]=0$.  


\begin{rema}
Changing notation, we could also have set $G=\overline{\rho(\pi_1(M))}$, the Zariski-closure of the image of the monodromy representation $\rho$ in $\PGL(\E)$, and $G_0=G\cap\PU(\E,h)$. Then by~\cite{Sim}, $G_0$ is the real Zariski closure of $\rho(\pi_1(M))$  and is a group of Hodge type. Moreover the subgroup of $G$ stabilizing the flag ${\mathbb F}^0\subset {\mathbb F}^{1}\subset\dots\subset {\mathbb F}^{k-1}\subset {\mathbb F}^k=\E$ is a parabolic subgroup $Q$ of $G$ and setting $V_0=G_0\cap Q$, 
we obtain another period domain $D=G_0/V_0$ in the flag manifold $G/Q$. Again, there is a period map $f$ from the universal cover of $M$ to $D$, which is $\rho$-equivariant, holomorphic and superhorizontal.
\end{rema}

\subsection{A possible strategy for proving Carlson-Toledo's conjecture} \label{conj}
In \cite{Re}, Reznikov treated a few special cases when $\G:=\pi_1(M)$ admits a rigid linear representation. More precisely, let $\rho:\G\fd\GL(m,\C)$ be a rigid, unbounded, reductive and irreducible representation. It follows from the general theory of Simpson~\cite{Sim} that such a representation necessarily is the monodromy of a complex variation of Hodge structure over $M$. As we just explained, if $G$, resp. $G_0$, is the Zariski closure, resp. real Zariski closure, of $\rho(\G)$ in $\GL(m,\C)$, then we obtain a period map $f$ from the universal cover $\tilde M$ of $M$ to a period domain $D$ of the form $G_0/V_0$. Reznikov assumed that $G_0$ is a Hermitian group and that in fact $D$ is the Hermitian symmetric space $X$ associated to $G_0$. Then the pull-back of the K\"ahler form $\o_X$ by $f$, which goes down to $M$ by equivariance, gives a non trivial cohomology class in $H^2(M,\R)$ (because $f^\star\o_X$ is semi-positive and $f$ is non constant). The fact that this class  comes from a class in $H^2(\G,\Z)$ follows from the fact that $\pi_2(X)=0$ (see below). 

When the period domain $D$ is not a bounded symmetric domain, for example if $G_0$ is not of Hermitian type, it might still be possible to produce a class in $H^2(\G,\R)$ from the geometry of the $\C$-VHS and the period map $f:\tilde M\fd D$. Indeed, there are $G_0$-homogeneous line bundles over $D$ (for example the canonical bundle $K_D$). Let ${\bf L}\fd D$ be such a line bundle and assume that it admits a $G_0$-invariant metric whose curvature form $\sqrt{-1}\,\Theta({\bf L})$ is positive on the superhorizontal distribution (this is the case for $K_D$ by Proposition~\ref{cano}). The line bundle $f^\star {\bf L}\fd \tilde M$ descends to a line bundle $L\fd M$ by equivariance and, as before, $\frac{\sqrt{-1}}{2\pi}f^\star \Theta({\bf L})$ induces a non trivial cohomology class $\eta\in H^2(M,\Z)$, which is just the first Chern class of $L$. Moreover, the image of $\eta$ in $H^2(M,\R)$ does not vanish because the semi-positivity of $\frac{\sqrt{-1}}{2\pi}f^\star \Theta({\bf L})$ implies that $\eta$ cannot be a torsion class. Again, it remains to prove that $\eta$ belongs to the image of the natural inclusion $\iota:H^2(\G,\Z)\hookrightarrow H^2(M,\Z)$.

In view of the exact sequence~\ref{long}: 
$$ 
\cdots \lo \pi_2(L^0) \lo \pi_2(M) \stackrel{\partial_L}{\lo}\pi_1(\cl^\star) = \z \lo
\pi_1(L^0) \lo \pi_1(M) \lo 1\, , 
$$
where $L^0$ is the $\C^\star$-bundle associated to $L$, it would be sufficient to show that the connecting morphism $\partial_L:\pi_2(M)\fd\Z$ is trivial for some line bundle $L\fd M$ induced by a homogeneous ${\bf L}\fd D$ as above. This is of course the case if $f_\star:\pi_2(M)\fd\pi_2(D)$ vanishes, in particular if $\pi_2(D)=0$, that is if $D=X$.






\subsection{Variations of Hodge structure with a rank 1 subbundle:
  proof of Theorem~\ref{cartoun}} \label{proofofcartoun}

In this section, we show that a slight modification of the strategy we
just described gives a proof of Theorem~\ref{cartoun}. 

We work here with the $\PU(p,q)$ period domain $D$ and map $f$ associated to the $\C$-VHS $E$ the representation $\rho$ comes from, and we keep the notation of Section~\ref{var}. 
Recall in particular that each $\theta_j:E^j\fd E^{j+1}\otimes\Omega^1_M$ deduced from the flat connection $\nabla$ on $E$ can be interpreted as a component of the differential of the period map.

By assumption, there exists $i$, $0<i<k$, such that $r_i={\rm rk}\,E^i=1$. First, we show that $\theta_{i-1}:E^{i-1}\fd E^i\otimes\Omega^1_M$ has (complex) rank 1 generically on $M$ (this kind of property was already observed by Siu, Carlson-Toledo and others). If it were not the case, $\theta_{i-1}$ would generically have rank greater than or equal to 2 ($\theta_{i-1}$ cannot vanish identically otherwise the variation of Hodge structure $E$ splits in two subvariations of Hodge structure and this is impossible since $\rho$ is assumed to be irreducible). 
This would imply that for a generic $x\in M$, writing $E^i_x=\C w$, there exist $v,v'\in E^{i-1}$ and $\g,\g'\in \Omega^1_{M,x}$ two linearly independent forms  such that $\theta_{i-1}(v)=w\otimes\g$ and  $\theta_{i-1}(v')=w\otimes\g'$. Since the forms $\g$ and $\g'$ are independent, $T_{M,x}=\Ker\g + \Ker\g'$. Now let $\zeta\in T_{M,x}$ and assume $\zeta\in\Ker\g$. Then, considering $\t_i$ as a morphism from $E^i\otimes T_M$ to $E^{i+1}$,  $\t_i(w\otimes\zeta)=\t_i((\t_{i-1}(v)\xi)\otimes\zeta)$ for some $\xi\in T_{M,x}$ such that $\g(\xi)=1$. By the commutation relation $[\t,\t]=0$, $\t_i(w\otimes\zeta)= \t_i((\t_{i-1}(v)\zeta)\otimes\xi)=\t_i((\g(\zeta)w)\otimes\xi)=0$. The same argument shows that $\t_i(w\otimes\zeta)=0$ for $\zeta\in\Ker\g'$, hence for all $\zeta$: $\t_i$ vanishes. In other words, if $\theta_{i-1}$ has (complex) rank greater than or equal to 2 then $\theta_{i}$ vanishes generically, hence everywhere. As above, this contradicts the irreducibility of $\rho$.   

Therefore, using again that $r_i=1$, the $\Hom(E^{i-1},E^i)$ component of $df$ also has complex rank 1 generically on $M$.


Consider the maximal parabolic subalgebra $\qg'\subset\gg$ associated to $\Pi_{\qg'}=\{\a_{r_0+\dots+r_{i-1}}\}$ (see Section~\ref{var}) and the corresponding parabolic subgroup $Q'\subset G$: we have $Q\subset Q'$. The projection map $\nu:D_u=G/Q\fd D'_u:=G/Q'$, which maps the flag $\{{F}^0\subset {F}^{1}\subset\dots\subset {F}^{k-1}\}$ to the $(r_0+\dots+r_{i-1})$-plane $F^{i-1}$, is holomorphic with respect to the natural complex structure on $D_u$ and $D'_u$. Its restriction to $D$, which we still denote by $\nu$, induces a holomorphic map $\nu:D\fd D'$ where $D'$ is the $G_0$-orbit of $eQ'$ in $D'_u$ (it is measurable open, see Section~\ref{mes}).

The map $g:=\nu\circ f:\tilde M\fd D'$ is holomorphic and since, by the very choice of $\nu$,  $d\nu$ kills the $\Hom(E^{j-1},E^j)$ component of $df$ for all $j\neq i$, $g$ has rank 1. As a consequence, the fibers of $g$ are pure dimensional, of codimension 1. We can thus apply a theorem of Kaup \cite[Theorem 6]{Gra} which states that there exists a smooth complex curve $N$ and holomorphic maps $\hat g:\tilde M\fd N$, $\psi:N\fd D'$, with $\hat g$ surjective, such that $g=\psi\circ \hat g$. We claim that $N$ cannot be isomorphic to $\P^1$. Indeed, $f$ is superhorizontal, hence the Levi form ${\cal L}(\varphi)$ is positive on each nonvanishing tangent vector in the image of the differential of $g$ (see Proposition \ref{canon} and the discussion following it). As a consequence, $\psi^\star\varphi$ is a subharmonic function on $N$ and its Levi form $\psi^\star{\cal L}(\varphi)$ is generically positive since $\psi$ is nonconstant. This is a contradiction if $N$ is compact.

Finally, we can apply the method described above, using $g^\star K_{D'}$ instead of $f^\star K_{D}$. The restriction of $g^\star K_{D'}$ to any 2-sphere lying in $\tilde M$ is clearly topologically trivial since $g$ factors through $N$ which is aspherical.

\section{Superhorizontal 2-spheres in period domains} \label{superhori}

Let $E$ be a $\C$-VHS over a compact K\"ahler manifold $M$, and $f:\tilde
M\fd D$ the corresponding period map. Since $f$ is superhorizontal, if
$\tilde\sigma:S^2\fd \tilde M$ is any (smooth) lift of some element
$[\sigma]\in\pi_2(M)$, $f\circ\tilde\sigma: S^2\fd D$ is a
superhorizontal 2-sphere in the period domain $D$. Of course, if such
superhorizontal spheres were homotopically trivial, then for every
homogeneous line bundle ${\bf L}\fd D$, the connecting homomorphism
$\partial_L: \pi_2(M)\fd\Z$ in the exact sequence~\ref{long} would be trivial and the non-vanishing of
$H^2(\pi_1(M),\R)$ would easily follow from the arguments of
Section~\ref{conj}. However, this is in general not true, and we will
show in this section that in most cases, there are many non-trivial
superhorizontal 2-spheres in the period domain $D$. In fact, when the hypotheses of 
Theorem~\ref{cartoun} are not met, and if $D$ is the $\PU(p,q)$ period domain associated to $E$, $\pi_2(D)$ is generated by
superhorizontal spheres, see Theorem~\ref{hori}.

\subsection{Construction of horizontal maps via Gromov $h$-principle}

Following~\cite{GromovPDR,GromovCC} and~\cite{Pansu}, we explain in
this section how the $h$-principle may be used to produce
horizontal maps. These considerations apply to a much wider setting
than superhorizontal maps in period domains, and the notation of
this section might differ from those used in the rest of the paper.     

Let $M$ be a real closed $k$-dimensional manifold, $N$ a real $n$-dimensional manifold, and let $H$ be a subbundle of the tangent bundle $T_N$ of $N$, called the horizontal distribution, ${\rm rk}(H)=h\geq k$. Let $\t$ be the projection $T_N\fd T_N/H$ seen as a 1-form on $N$ with values in $T_N/H$ and let $d\t:H\times H\fd T_N/H$ be defined on horizontal vector fields by $d\t(\xi,\eta)=[\xi,\eta]\mod H$ ($d\t$ is in fact a 2-form on $H$ with values in $T_N/H$).

Assume that $f$ is an $H$-horizontal immersion of $M$ into $N$, namely an immersion such that for all $x\in M$, $d_xf$ maps $T_{M,x}$ in $H_{f(x)}\subset T_{N,f(x)}$. The differential of $f$ gives rise to a map $F:M\fd {\rm Gr}_{k}(H)$, where ${\rm Gr}_{k}(H)\fd N$ is the bundle of Grassmannians of $k$-planes in $H$ over $N$.  The horizontality of $f$ is equivalent to the vanishing $f^\star\t=0$. Moreover, if $f$ is $C^2$, then $f^\star d\t=0$ also holds. This means that $df$ sends the tangent spaces $T_{M,x}$ to {\em isotropic} subspaces of $H$, namely subspaces $S$ of $H_{f(x)}$ on which $d\t$ vanishes. Therefore, a necessary condition for a $C^2$ horizontal immersion to exist is the existence of a map $F$ from $M$ to the bundle of Grassmannian of horizontal isotropic $k$-planes over $N$.

We shall indicate the different steps of the construction of a {\em folded} horizontal immersion $f:M\fd N$ from the data of $F:M\fd {\rm Gr}_{k}^{reg,iso}(H)$, where  ${\rm Gr}_{k}^{reg,iso}(H)\fd N$ is the bundle of {\em regular} horizontal isotropic $k$-planes over $N$. The regularity assumption will be explained in the next paragraphs. The meaning of a folded immersion will be given later.  

\bigskip

First, one wants to construct germs of horizontal immersions tangent to given horizontal subspaces of $T_N$ at given points. One therefore wants to locally solve the horizontality equation ${\mathcal H}(f):=f^\star\t=0$. Since this is a local problem, we assume that $M$ is an open set in $\R^k$, $0\in M$, and that the quotient bundle $T_N/H$ is trivial on $N$, so that $\t$ and $d\t$ can be seen as globally defined $\R^{n-h}$-valued forms on $N$.  
Here, ${\mathcal H}$ is understood as a differential operator from the space of sections of the trivial bundle $M\times N\fd M$ (so that sections of this bundle are just maps $M\fd N$) to the space of $\R^{n-h}$-valued 1-forms on $M$. Spaces of sections are endowed with the {\em fine} topology: a basis of the fine $C^0$-topology on the space of continuous sections $C^0(M,E)$ of a bundle $E$ over $M$ is given by the sets $C^0(M,{\mathcal V})$ of continuous ${\mathcal V}$-valued sections of $E$, for all open sets ${\mathcal V}\subset E$. The fine $C^r$-topology on $C^r(M,E)$ is induced from the fine $C^0$-topology on the space of sections of the $r$-jets bundle $E^{(r)}\fd M$ via the natural inclusion $C^r(M,E)\subset C^0(M,E^{(r)})$.

\medskip

The idea is to right invert the linearization $D_f{\mathcal H}$ of the operator ${\mathcal H}$ at $f$. If $\xi$ is a section of $f^\star T_N$, then 
$$
D_f{\mathcal H}(\xi)=d(\iota_\xi\t)+f^\star(\iota_\xi d\t).
$$  
Therefore, if $\xi$ is horizontal (i.e. $\xi$ is a section of $f^\star H$),
$D_f{\mathcal H}(\xi)=f^\star(\iota_\xi d\t)$. The linearization $D_f{\mathcal H}$ is seen as a bundle map
from $f^\star H$ to the space $M\times\Hom(\R^k,\R^{n-h})$ of $\R^{n-h}$-valued 1-forms on $M$. That is where the regularity assumption comes from: a $k$-dimensional subspace $S$ of $T_{N,y}$ is said to be {\em regular} if the linear map 
$$
\begin{array}{rcl}
H_{y} & \fd & \Hom(S,\R^{n-h})\\
v & \longmapsto & \iota_v d\t
\end{array}
$$ 
is surjective (if $S\subset T_{N,y}$ is not horizontal, i.e. not contained in $H_y$, this definition depends on the choice of $\t$ as an $\R^{n-h}$-valued 1-form on $N$). 

The set of regular $k$-dimensional subspaces is open in the bundle of Grassmannians of $k$-planes in $T_N$, 
and so is the set ${\rm RegIm}(M,N)$ of regular immersions $M\fd N$ in the space of maps $M\fd N$.

It turns out that the regularity assumption is enough for being able to solve
the horizontality equation up to infinite order, and then to right invert the
linearization of ${\mathcal H}$. Namely:

\begin{prop}{\rm (\cite[Prop. 29]{Pansu})}.\label{allorders} Let $S\subset H_y$ be a $k$-dimensional
regular isotropic subspace. Then there exists a germ of immersion $f:(\R^k,0)\fd
(N,y)$ which is tangent to $S$ at the origin: $d_0f(\R^k)=S$, and which
satisfies the horizontality equation to all orders: ${\mathcal H}(f)(x)=o(|x|^m)$ for all
$m\in\N$. 
\end{prop}

\begin{prop}{\rm (See~\cite[Prop. 31]{Pansu})}.\label{rightinvertible} Assume $f:M\subset\R^k\fd N$ is a regular immersion (i.e. $d_xf(\R^k)$ is a
regular $k$-dimensional subspace of $T_{N,f(x)}$ for all $x\in M$), then the map 
$$
\begin{array}{rcl}
	\Gamma(f^\star H) & \fd & \Gamma(M\times\Hom(\R^k,\R^{n-h}))\\
	\xi & \longmapsto & f^\star(\iota_\xi d\t) 
\end{array}
$$
is right invertible.
\end{prop}

\begin{demo}
The set of regular $k$-planes in the bundle of Grassmannians
${\rm Gr}_k(T_N)$ is open as we mentionned. Hence the set ${\rm RegInj}(T_M,T_N)$ of triples $(x,y,L)$ with $x\in M$, $y\in N$ and $L$ an injective linear map from $\R^k$ to $T_{N,y}$ with regular image, is open in the bundle $\Hom(T_M,T_N)$ over $M\times N$. 

Given $(x,y,L)\in {\rm RegInj}(T_M,T_N)$, the set of right inverses of the linear map 
$$
\begin{array}{rcl}
H_y & \fd & \Hom(\R^k,\R^{n-h}) \\
v & \longmapsto & (d\t)_y(v,L(.)) \\
\end{array}
$$ 
is an affine space of dimension $k(n-h)(h-k(n-h))$. Therefore the collection of these spaces over ${\rm RegInj}(T_M,T_N)$ is a bundle, with contractible fibers, and hence admits a section $(x,y,L)\mapsto {\rm rightinv}(x,y,L)$. 

Now, given an $\R^{n-h}$-valued 1-form $\beta$ on $M$, i.e. $\beta(x)\in\Hom(\R^k,\R^{n-h})$ for all $x\in M$, define an horizontal vector field along $f$ by $\xi_f(\beta)=(x,{\rm rightinv}(x,f(x),d_xf)(\beta(x)))$. $\xi_f$ is a right inverse of $D_f{\mathcal H}$. 
\end{demo}

\medskip

Now, a theorem of Gromov in the spirit of the inverse function theorem of Nash shows that if the linearization $D_f{\mathcal H}$ of the operator ${\mathcal H}$ can be right inverted for all $f\in{\rm RegIm}(M,N)$, then locally, the operator ${\mathcal H}$ itself can be right inverted: for any regular immersion $f$, there exists a right inverse ${\mathcal H}_f^{-1}$ of ${\mathcal H}$, defined on a $C^s$-neighborhood ${\mathcal V}_f$ of ${\mathcal H}(f)$, such that ${\mathcal H}_f^{-1}({\mathcal H}(f))=f$. Moreover, ${\mathcal H}_f^{-1}$ depends smoothly on parameters. See~\cite[p. 114-120]{GromovPDR} for a precise statement of the theorem and its corollaries in the form that is useful here.

This theorem and the two propositions prove that the Cauchy problem for horizontal immersions, with initial data a point $y\in N$ and a regular horizontal isotropic $k$-plane $S$ in $T_{N,y}$, has a local solution: there exists a germ of horizontal immersion $f:(\R^k,0)\fd (N,y)$ which is tangent to $S$ at the originAbelian. 
Indeed, Proposition~\ref{allorders} gives a germ of regular immersion $f_0:(\R^k,0)\fd(N,y)$ such that ${\mathcal H}(f_0)=o(|x|^m)$. Proposition~\ref{rightinvertible} and Gromov's inverse function theorem~\cite[Main Theorem (3), p. 117]{GromovPDR} shows that there exists a right inverse ${\mathcal H}_{f_0}^{-1}$ of ${\mathcal H}$ defined in a neighborhood of ${\mathcal H}(f_0)$. Taking a smaller neighborhood if necessary, we may assume that the identically zero 1-form belongs to this neighborhood. Applying ${\mathcal H}_{f_0}^{-1}$  gives the desired germ $f$. 

In fact, we can do better: the germ $f_0$ can be deformed to a true horizontal germ~\cite[Corollary (D), p. 119]{GromovPDR}. Indeed, again working on a smaller neighborhood of $0\in\R^k$ if necessary, we may assume that $(1-t){\mathcal H}({f_0})$ belongs to the given neighborhood ${\mathcal V}_{f_0}$ of ${\mathcal H}(f_0)$ for all $t\in[0,1]$ and we have the deformation $f_t:={\mathcal H}_{f_0}^{-1}((1-t){\mathcal H}({f_0}))$ from $f_0$ to $f_1$ which is horizontal. 

We will say that regular horizontal immersions satisfy the {\em local $h$-principle}.

\medskip 

The right inverse ${\mathcal H}_f^{-1}$ has moreover the property of being a local operator~\cite[Main Theorem (5), p. 118]{GromovPDR}: for a given Riemannian metric on $M$ and $r>0$, we may assume that for all regular immersions $f:M\fd N$, all 1-forms $\beta\in {\mathcal V}_f$ and all $x\in M$, the value of ${\mathcal H}_f^{-1}(\beta)$ at $x$ depends only on the values of $f$ and $\beta$ in a ball of radius $r$ around $x$ (of course ${\mathcal V}_f$ depends on the metric and on $r$). This in turns implies an important fact about regular horizontal immersions $M\fd N$: they are {\em microflexible}~\cite[p. 120]{GromovPDR}. By definition, this means that given two compact sets $K'\subset K\subset M$, a continuous family of regular horizontal immersions $\{f_a\}_{a\in A}$ defined near $K$ and indexed by a polyhedron $A$, and a family $\{f_{a,t}\}_{a\in A}$ of regular horizontal immersions, defined near $K'$, continuous in $a\in A$ and $t\in[0,1]$, such that $f_{a,0}=f_{a}$ near $K'$ for all $a\in A$, there exist $\varepsilon>0$ and  a continuous family $\{\tilde f_{a,t}\}_{a\in A}$ of regular horizontal maps, defined near $K$ and for $t\in[0,\varepsilon]$, such that for all $a\in A$, $\tilde f_{a,0}=f_a$ near $K$ and $\tilde f_{a,t}=f_{a,t}$ near $K'$. 

Indeed, for $t$ small, we can find regular immersions $F_{a,t}$ from a neighborhood of $K$ in $M$ to $N$ such that $F_{a,t}=f_{a,t}$ near $K'$, and then, for $t$ small, we may define $\tilde f_{a,t}={\mathcal H}_{F_{a,t}}^{-1}(0)$. This gives a family of regular horizontal immersions defined near $K$ which agree with $f_{a,t}$ near $K'$ by locality. 

\bigskip

By~\cite[p. 262]{GromovCC}, the local $h$-principle and the microflexibility imply that given our $k$-dimensional closed manifold $M$ and a continuous map $F$ from $M$ to the bundle ${\rm Gr}_{k}^{reg,iso}(H)\stackrel{\pi}{\fd} N$ of Grassmannians of regular $H$-horizontal isotropic $k$-planes over $N$, there is a polyhedron $W$, homotopy equivalent to $M$, and a folded horizontal immersion $f^\vee:W\fd N$ such that $f^\vee\circ e$ approximates $\pi\circ F$ (for the topology of uniform convergence and some homotopy equivalence $e:M\fd W$). Here by a folded horizontal immersion $f^\vee:W\fd N$ we mean that $f^\vee$ is continuous and is a smooth horizontal immersion on each simplex of $W$ (up to the boundary). 

The construction of such a map is explained in~\cite{Pansu} when $M=S^1$. In the proof of the following proposition, we will give a detailed construction for $M=S^2$ since this is the case we are ultimately interested in. For the general case, see Remark~\ref{remarque1}.

\begin{prop}\label{folded}
Let $f_0:S^2\fd N$ be a continuous map and assume it admits a continuous lift $F$ to the bundle ${\rm Gr}_{2}^{reg,iso}(H)\stackrel{\pi}{\fd} N$ of Grassmannians of regular $H$-horizontal isotropic $2$-planes over $N$, such that $\pi\circ F$ is homotopic to $f_0$. 
Then there exists a smooth horizontal map $f:S^2\fd N$ homotopic to $f_0$.
\end{prop}

\begin{demo}
We begin with the construction of a polyhedron $W$, homeomorphic to $S^2$, and a folded horizontal immersion $f^\vee:W\fd N$, which ``approximates'' the map $f_0$. 
From our discussion so far, we know that there exists a continuous choice of germs of horizontal immersions $f_x$, $x\in S^2$,  defined on a neighborhood of $0\in\R^2$, such that $f_x(0)=\pi\circ F(0)\in N$ and $d_0f_x(\R^2)=F(x)$. We may assume that all these germs are defined in the Euclidean ball ${\bf B}$ of radius 2 around 0. Consider a fixed equilateral triangle ${\bf T}$ inscribed in the circle of radius 1 in ${\bf B}$, let ${\bf v}_r$, ${\bf v}_g$ and ${\bf v}_b$ be its vertices, and ${\bf e}_{rg}$, ${\bf e}_{gb}$, ${\bf e}_{br}$ its edges. Using microflexibility, we are going to construct families of horizontal immersions defined near ${\bf T}$ in ${\bf B}$, with prescribed boundary values on $\partial{\bf T}$, from which we shall select maps which we shall then glue together along a pattern given by a suitably chosen triangulation of $S^2$. 

All the forthcoming metric notions on the sphere $S^2$ are to be understood with respect to the round metric. 
\medskip

Let ${\bf U}_r$, ${\bf U}_g$, ${\bf U}_b$ be small disjoint discs centered at ${\bf v}_r$, ${\bf v}_g$, ${\bf v}_b$, and ${\bf U}_{rg}$, ${\bf U}_{gb}$, ${\bf U}_{br}$ be small cylindrical neighborhoods of ${\bf e}_{rg}$, ${\bf e}_{gb}$, ${\bf e}_{br}$ in ${{\bf B}}$ such that ${\overline{\bf U}_{rg}\cap \overline{\bf U}_{gb}}\subset {\bf U}_g$, ${\overline{\bf U}_{gb}\cap \overline{\bf U}_{br}}\subset {\bf U}_b$ and ${\overline{\bf U}_{br}\cap \overline{\bf U}_{rg}}\subset {\bf U}_r$.

Let $x\in S^2$ and $\xi$ be a unit tangent vector to $S^2$ at $x$. For $t\in[0,1]$, define an horizontal immersion $f_{x,\xi,t}$ on a neighborhood of $\overline{\bf U}_r\cup \overline{\bf U}_g$ by 
$$
f^{rg}_{x,\xi,t}=\left\{
\begin{array}{l}
 f_{c_{x,\xi}(-t)} \mbox{ near } \overline{\bf U}_r\\
 f_{c_{x,\xi}(t)} \mbox{ near } \overline{\bf U}_g\\
\end{array}
\right.
 $$
where $c_{x,\xi}$ is the unit speed geodesic going through $x$ in the direction of $\xi$. 

The maps $f^{rb}_{x,\xi,t}$, $t\in[0,1]$, give a deformation of the horizontal map $f_x:{{\bf B}}\fd N$ near $\overline{\bf U}_r\cup \overline{\bf U}_g$. By microflexibility, there exists an $\varepsilon_0 >0$, uniform in $x$ and $\xi$, and a deformation of $f_x$ by horizontal immersions $\tilde f^{rg}_{x,\xi,t}:{{\bf B}}\fd N$, $t\in[0,\varepsilon_0]$, such that 
$\tilde f^{rg}_{x,\xi,t}=f^{rg}_{x,\xi,t}$ near $\overline{\bf U}_r\cup \overline{\bf U}_g$. 

In the same way we define horizontal maps $\tilde f^{gb}_{x,\xi,t}$ and $\tilde f^{br}_{x,\xi,t}$, for $t\in[0,\varepsilon_0]$.

Let $\xi_r$, $\xi_g$ and $\xi_b$ be unit tangent vectors at $x\in S^2$. Let $m_{rg}(t)$, $m_{gb}(t)$ and $m_{br}(t)$ be the midpoints of the geodesic segments $[c_{x,\xi_r}(t),c_{x,\xi_g}(t)]$, $[c_{x,\xi_g}(t),c_{x,\xi_b}(t)]$ and  $[c_{x,\xi_b}(t),c_{x,\xi_r}(t)]$. Let also $\xi_{rg}(t)$, $\xi_{gb}(t)$, $\xi_{br}(t)$ and $s_{rg}(t)$, $s_{gb}(t)$, $s_{br}(t)$ be such that 
$$
\left\{
\begin{array}{l}
c_{m_{rg}(t),\xi_{rg}(t)}(s_{rg}(t))=c_{x,\xi_g}(t)\\ 
c_{m_{gb}(t),\xi_{gb}(t)}(s_{gb}(t))=c_{x,\xi_b}(t)\\ 
c_{m_{br}(t),\xi_{br}(t)}(s_{br}(t))=c_{x,\xi_r}(t)
\end{array}
\right.
$$
Everything is made so that the maps $f^{rgb}_{x,\xi_r,\xi_g,\xi_b,t}$ given for $t\in[0,\varepsilon_1]$ for some $0<\varepsilon_1\leq \varepsilon_0$ by
$$
f^{rgb}_{x,\xi_r,\xi_g,\xi_b,t}=\left\{
\begin{array}{l}
 \tilde f^{rg}_{m_{rg}(t),\xi_{rg}(t),s_{rg}(t)} \mbox{ near } \overline{\bf U}_{rg}\\
 \tilde f^{gb}_{m_{gb}(t),\xi_{gb}(t),s_{gb}(t)} \mbox{ near } \overline{\bf U}_{gb}\\
\tilde f^{br}_{m_{br}(t),\xi_{br}(t),s_{br}(t)} \mbox{ near } \overline{\bf U}_{br}\\
\end{array}
\right.
$$
are well defined on $\overline{\bf U}_{rg}\cup \overline{\bf U}_{gb}\cup \overline{\bf U}_{br}$. 
Indeed, for example on $\overline{\bf U}_{rg}\cap \overline{\bf U}_{gb}\subset {\bf U}_g$, we have 
$$
\begin{array}{rcl}
\tilde f^{rg}_{m_{rg}(t),\xi_{rg}(t),s_{rg}(t)} & = & f_{c_{m_{rg}(t),\xi_{rg}(t)}(s_{rg}(t))}\\
& = & f_{c_{x,\xi_g}(t)}\\
& = & f_{c_{m_{gb}(t),\xi_{gb}(t)}(-s_{gb}(t))}\\
& = & \tilde f^{gb}_{m_{gb}(t),\xi_{gb}(t),s_{gb}(t)}   
\end{array}
$$

Hence these maps give a deformation of $f_x:{{\bf B}}\fd N$ near $\overline{\bf U}_{rg}\cup \overline{\bf U}_{gb}\cup \overline{\bf U}_{br}$ and by microflexibility, there exist $\varepsilon>0$, uniform in $x$, $\xi_r$, $\xi_g$, $\xi_b$, and a deformation of $f_x$ by horizontal immersions  $\tilde f^{rgb}_{x,\xi_r,\xi_g,\xi_b,t}:{{\bf B}}\fd N$ defined for $t\in[0,\varepsilon]$ such that 
$\tilde f^{rgb}_{x,\xi_r,\xi_g,\xi_b,t}=f^{rgb}_{x,\xi_r,\xi_g,\xi_b,t}$ near $\overline{\bf U}_{rg}\cup \overline{\bf U}_{gb}\cup \overline{\bf U}_{br}$.

\medskip

Now consider a triangulation ${\mathcal T}$ of $S^2$, and assume that it is geodesic, that each triangle $T$ of $\mathcal T$ is contained in a convex ball of radius $<\varepsilon/p$ for some large number $p$, and that its angles are not too obtuse (the center of the circumcircle of $T$ should also be contained in the ball of radius $<\varepsilon/p$ containing $T$). Assume moreover 
that $\mathcal T$ is even, namely that there is an even number of triangles around each vertex of $\mathcal T$. Triangulations satisfying our assumptions exist for any value of $p$, as is easily verified. The fact that the triangulation is even 
implies that $\mathcal T$ is 3-colorable, i.e. that we can color its vertices with 3 colors, say red, green and blue, so that adjacent vertices have different colors. Abelian

Choose a 3-coloring of the vertices of the triangulation ${\mathcal T}$ and,
for each triangle $T$ of $\mathcal T$, let $x^T$ be the center of the circumcircle of $T$  and let $\xi^T_r$ (resp. $\xi^T_g$, $\xi^T_b$) be the unit tangent vector at $x^T$ in the direction of the red (resp. green, blue)  vertex of $T$. 
Construct the horizontal maps $f_{T,t}=\tilde f^{rgb}_{x^T,\xi^T_r,\xi^T_g,\xi^T_b,t}$ as before. Let $t^T$ be such that $m_{rg}(t^T)$ (resp. $m_{gb}(t^T)$, $m_{br}(t^T)$) is the midpoint of the red-green (resp. green-blue, blue-red) edge of $T$ and set $f_T=f_{T,t^T}$. The maps $f_T$ are well-defined because we assumed that the triangles are very small.
Moreover, if two triangles $T$ and $T'$ are adjacent, say along their red-green edges, then the maps $f_T$ and $f_{T'}$ agree on ${\bf U}_{rg}\subset{ {\bf B}}$ because they depend only on the red and green vertices they share, and on the midpoint of the red-green edge connecting them. 

Consider one copy ${\bf T}_T$ of ${\bf T}$ for each triangle $T$ of the triangulation $\mathcal T$, and identify the red-green (resp. green-blue, blue-red) edge of ${\bf T}_T$ with the red-green (resp. green-blue, blue-red) edge of ${\bf T}_{T'}$ if $T$ and $T'$ intersect along their red-green (resp. green-blue, blue-red) edges. The resulting polyhedron $W$ is homeomorphic to the sphere $S^2$. Moreover, the map $f^\vee:W\fd N$ given by $f^\vee=f_T$ on ${\bf T}_T$ is well-defined, and it is the folded horizontal immersion we wanted.     

\medskip

Neither the polyhedron $W$ nor the map $f^\vee$ are smooth, but we can easily build a smooth horizontal map $f:S^2\fd N$ out of the $f_T$'s (this $f$ will not be an immersion). 

For all vertices $v$ of the triangulation ${\mathcal T}$ of $S^2$, let $B_v$ be the ball of radius $r$ around $v$, $r$ being such that these balls are pairwise disjoint. For all edges $e:[0,L_e]\fd S^2$ of ${\mathcal T}$ parameterized by arc length, let $\eta_e$ be a unit normal vector field along $e$ and let $U_e$ be the open ``tubular'' set $\{c_{e(t),\eta_e(t)}(s)\ |\ t\in(r/2,L_e-r/2),\ s\in(-r',r')\}$ for some $r'<r/2$, chosen so that these sets are pairwise disjoint. 

Given an edge $e$ of ${\mathcal T}$, we map the set $B_{e(0)}\cup U_e\cup B_{e(L_e)}$ to the corresponding edge of the model triangle ${\bf T}$. Say $e(0)$ is a red vertex and $e(L_e)$ is a green one. Define a map $\phi_e$ from the edge $e$ to the red-green edge ${\bf e}_{rg}:[0,\sqrt{3}]\fd \R^2$ of ${\bf T}$ as follows
$$
\phi_e:\ e(t)\longmapsto\left\{
\begin{array}{lr}
{\bf e}_{rg}(0)={\bf v}_r & \mbox{ if }t\in[0,r]\\
{\bf e}_{rg}(\chi_e(t)) & \mbox{ if }t\in[r,L_e-r]\\
{\bf e}_{rg}(\sqrt{3})={\bf v}_g & \mbox{ if }t\in[L_e-r,L_e]
\end{array}
\right.
$$ 
where $\chi_e : [r,L_e-r]\fd [0,\sqrt{3}]$ is a smooth, onto, strictly increasing function all of whose derivatives vanish at $r$ and $L_e-r$. Now let $\Phi_e$ by the map from $B_{e(0)}\cup U_e\cup B_{e(L_e)}$ to ${\bf e}_{rg}$ given by:
 $$
\Phi_e(x)=\left\{
\begin{array}{lr}
{\bf v}_r & \mbox{ if }x\in B_{e(0)}\\
{\phi_e(e(t))}\in{\bf e}_{rg} & \mbox{ if }x=c_{e(t),\eta_e(t)}(s)\in U_e\\
{\bf v}_g & \mbox{ if }x\in B_{e(L_e)}
\end{array}
\right.
$$ 
The $\Phi_e$'s are smooth and they agree on $B_v$ whenever two edges are adjacent at a vertex $v$. Hence, if $V$ is the neighborhood of the $1$-skeleton of ${\mathcal T}$ given by the union of the $B_v$'s and the $U_e$'s, we obtain a smooth map $\Phi_V$ from $V$ to the boundary $\partial{\bf T}$ of ${\bf T}\subset\R^2$. Multiplying it by a smooth map $S^2\fd [0,1]$ which is identically 0 outside $V$ and identically 1 inside a smaller neighborhood of the 1-skeleton, we get a smooth map $\Phi:S^2\fd {\bf T}$, which sends the vertices and edges of the triangulation ${\mathcal T}$ to the vertices and edges of ${\bf T}$ according to the chosen $3$-coloring of the triangulation (note that $\Phi$ sends $S^2\backslash V$ to the center $0\in\R^2$ of $\bf T$). Moreover, by construction, all normal derivatives of $\Phi$ at the edges of ${\mathcal T}$ vanish (and $\Phi$ is constant near the vertices). This implies that setting $f_{|T}=f_T\circ \Phi$ defines a smooth horizontal map $f:S^2\fd N$.        

\medskip

If the triangulation ${\mathcal T}$ is chosen sufficiently fine, namely if $p$ is sufficiently big, $f$ is very close to the map $\pi\circ F$ and therefore homotopic to it, hence to $f_0$.
\end{demo}

\begin{rema}\label{remarque1}
 If $M$ is any closed $k$-dimensional manifold, the arguments given here can be adapted to prove the existence of a folded horizontal immersion (or a smooth horizontal map) $M\fd N$. It is indeed enough to have a triangulation of $M$ with an  identification of its simplices with some model simplex, which will be used to glue horizontal immersions defined on copies of this model simplex, as we just did. This identification is given by a proper $(k+1)$-coloring of the vertices of the triangulation. Such triangulations, as fine as needed, always exist, since the barycentric subdivision of any triangulation of $M$ is $(k+1)$-colorable.   
\end{rema}

\begin{rema}\label{remarque2}
 If instead of a map $F:M\fd {\rm Gr}^{reg,iso}_{k}(H)$, we start with a map $F:M\fd {\rm Gr}^{reg,iso}_{k,l}(H)$ in the bundle of flags of regular $k$-planes in regular horizontal isotropic $l$-planes, $l>k$, then the so-called micro-extension trick can be used to get a smooth (non folded) horizontal immersion $M\fd N$, see~\cite[p. 258]{GromovCC} (and~\cite[Theorem 4.6.1]{EM} for a similar statement). 
\end{rema}

\subsection{Application to superhorizontal spheres in period domains}

Using the results of the previous section, we shall now prove that in general period domains have plenty of non-trivial superhorizontal 2-spheres. 

\subsubsection{A model case}

We begin by studying the period domain $D= \PU(2,n)/{\rm P}(\U(1)\times \U(n)\times \U(1))$.
It is an open set in the flag manifold $D_u=\PU(n+2)/{\rm P}(\U(1)\times
\U(n)\times \U(1))={\rm PGL}(n+2,\C)/Q$, where $Q$ is the parabolic subgroup of
${\rm PGL}(n+2,\C)$ stabilizing the flag $\{\la e_1\ra\subset\la e_1,\ldots,e_{n+1}\ra\subset
\C^{n+2}\}$ of $\C^{n+2}$. 
The complex dimension of $D$ is $2n+1$ and the superhorizontal distribution $H_D$ of $D$ is a holomorphic subbundle of $T_D$ of complex rank $2n$. 


We want to construct {\em superhorizontal} maps $S^2\fd D$, homotopic to given elements of the second homotopy group $\pi_2(D)$ of $D$. Note that $\pi_2(D)\simeq \Z$, as follows from the discussion in~\ref{bur}. A map $f:S^2\fd D$ is said to be superhorizontal simply if it is $H_D$-horizontal, namely if for all $x\in S^2$ and $\xi\in T_{S^2,x}$, $f_\star\xi\in H_{D,f(x)}$. 


\begin{prop}\label{petit} Let $D$ be the period domain $\PU(2,n)/{\rm P}(\U(1)\times \U(n)\times \U(1))$.
 For $n\geq 2$, all elements of $\pi_2(D)$ can be represented by smooth superhorizontal spheres $S^2\fd D$. For $n\geq 3$, they can be represented by smooth immersed superhorizontal spheres. 
\end{prop}

\begin{demo}
By Proposition~\ref{folded}, all we have to prove is that the elements of $\pi_2(D)$ can be lifted to $\pi_2({\rm Gr}^{reg,iso}_{2}(H_D))$, where ${\rm Gr}^{reg,iso}_{2}(H_D)$ is the bundle of Grassmannians of regular isotropic superhorizontal 2-planes over $D$. Note that $H_D\fd D$ is a holomorphic complex vector bundle, but we consider {\em real} $2$-planes, namely  real spans of couples of $\R$-linearly independent vectors of $H_D$. Regularity and isotropy are defined as before with respect to the 2-form $d\t:H_D\times H_D\fd T_D/H_D$. At the point $eQ\in D\subset D_u$, $H_D$ can be identified with $\C^{2n}=\C^{n}\times\C^n$ and $d\t$ restricts to the complex symplectic form 
$$
\begin{array}{rcl}
 q:(\C^{n}\times\C^n)\times (\C^{n}\times\C^n) & \fd & \C \\
((v_1,v_2),(w_1,w_2)) & \longmapsto & \trans v_1 w_2-\trans v_2 w_1
\end{array}
$$

Since the superhorizontal 2-planes in $T_D$ are candidates for being the images of the tangent
planes of $S^2$ by the differential of a map $f:S^2\fd D$, and since we have a natural orientation on $T_{S^2}$, we will in fact consider the bundle ${\rm Gr}^{reg,iso}_{2,or}(H_D)$ of regular isotropic superhorizontal {\em oriented} 2-planes.
If we can lift elements of $\pi_2(D)$ to elements of $\pi_2({\rm Gr}^{reg,iso}_{2,or}(H_D))$, then using the natural projection 
${\rm Gr}^{reg,iso}_{2,or}(H_D)\fd {\rm Gr}^{reg,iso}_{2}(H_D)$, we are done. 

\medskip

The fiber of the fibration ${\rm Gr}_{2,or}^{reg, iso}(H_D)\fd D$ at $eQ$ is ${\rm Gr}_{2,or}^{reg, iso}$, the Grassmannian of oriented (real) $2$-planes $S$ in $\C^{2n}$, isotropic for the complex symplectic form $q$ and regular in the sense that
$$
\begin{array}{rcl}
\C^{2n} & \fd & \Hom_\R(S,\C)\\  v & \longmapsto & (\iota_v{q})_{|S}
\end{array}
$$
is surjective. It is easy to check that a real $2$-plane $S$ in $\C^{2n}$ is regular if and only if it isn't a complex line, namely if and only if it is spanned by two $\C$-linearly independent vectors. 

The long exact sequence of homotopy groups associated to the fibration is
$$
\cdots \fd \pi_2({\rm Gr}_{2,or}^{reg, iso}(H_D)) \fd \pi_2(D) \fd \pi_1({\rm Gr}_{2,or}^{reg,iso}) \fd \cdots
$$
Therefore the elements of $\pi_2(D)$ which come from elements of $\pi_2({\rm Gr}_{2,or}^{reg, iso}(H_D))$ are those which are in the kernel of the map $\pi_2(D) \fd\pi_1({\rm Gr}_{2,or}^{reg, iso})$.  

\medskip

For $k\leq n$, the $q$-symplectic group ${\rm Sp}(n,\C)$ acts transitively on the set of $k$-tuples of $\C$-linearly
independent vectors spanning a $k$-dimensional complex isotropic subspace in $\C^{2n}$, and the stabilizer $Q_k$ of such a $k$-tuple is topologically ${\rm Sp}(n-k,\C)\times\C^{2k(n-k)}\times\C^{{k(k+1)}/{2}}$, as a simple computation shows. 

Therefore for $n\geq 2$ the stabilizer in ${\rm Sp}(n,\C)$ of an oriented regular isotropic $2$-plane is isomorphic to $\GL^+(2,\R)\times Q_2$. Hence ${\rm Gr}_{2,or}^{reg,iso}\simeq{\rm Sp}(n,\C)/(\GL^+(2,\R)\times Q_2)$. Since ${\rm Sp}(n,\C)$ is path-connected and simply connected, $\pi_1({\rm Gr}_{2,or}^{reg,iso})\simeq \pi_0(\GL^+(2,\R)\times Q_2)$. Now $\GL^+(2,\R)\times Q_2$ is path-connected, hence ${\rm Gr}_{2,or}^{reg,iso}$ is simply connected and all the elements of $\pi_2(D)$ can be lifted to elements of $\pi_2({\rm Gr}_{2,or}^{reg,iso}(H_D))$.

\medskip

If $n\geq 3$, we may also look at flags of oriented $2$-planes in oriented superhorizontal regular isotropic $3$-planes:
${\rm Gr}_{2,3,or}^{reg,iso}\simeq{\rm Sp}(n,\C)/((\GL^+(2,\R)\times\R_+^\star\times\R^2)\times Q_3)$ and find that all
elements of $\pi_2(D)$ can be lifted to elements of
$\pi_2({\rm Gr}_{2,3,or}^{reg,iso}(H_D))$. This implies that for $n\geq 3$, all the elements of $\pi_2(D)$ can be represented by smooth immersed superhorizontal spheres, see Remark~\ref{remarque2}.
\end{demo}

\begin{rema}\label{n=1}
 If $n=1$, one sees from the proof above that there are no regular isotropic horizontal real $2$-planes (indeed in this case a real 2-plane $S$ is isotropic if and only if it is a complex line if and only if it is not regular), hence the $h$-principle doesn't apply (at least in its present form).
\end{rema}

\subsubsection{Second homotopy groups of $\PU(p,q)$ period domains and flag manifolds} \label{bur}
Keeping the notation of Section~\ref{var}, let $D$ be the $\PU(p,q)$ period
domain associated to some variation of Hodge structure
$E=\bigoplus_{i=0}^k E^i$ over a compact K\"ahler manifold $M$ and let $D_u=\PU(p+q)/{\rm P}(\Pi_{i=0}^k{\rm U}(r_i))$ be its dual.  

We first sum up some results of~\cite[Chap. 4.D]{BR}.
Let $\gg^\a\subset\ng$ be a root-space (i.e. $\a\in\Phi$ and $-\a\not\in\Phi$). Then $\gg^\a\oplus\gg^{-\a}\oplus[\gg^\a,\gg^{-\a}]$ is a Lie algebra isomorphic to $\sg\lg(2,\C)$ which is the complexification of a subalgebra $\sg_\a$ of $\gg_u$ isomorphic to $\sg\ug(2)$. Let $S_\a$ be the corresponding copy of $\PU(2)$ and $U_\a\subset S_\a$ the connected subgroup with Lie algebra $\ug_\a=[\gg^\a,\gg^{-\a}]\cap\gg_0\subset\vg_0$. Then $U_\a$ is a copy of $U(1)$ so that $S_\a/U_\a$ is a 2-sphere. The homomorphism $i_\a:S_\a\fd G_u$ corresponding to the inclusion $\sg_\a\subset\gg_u$ defines an $i_\a$-equivariant immersion $\varphi_\a:\P^1\fd D_u$ (since $i_\a(U_\a)\subset V_0$) which is moreover totally geodesic and holomorphic.

The classes $[\Sigma_\a]$ of the spheres $\Sigma_\a:=\varphi_\a(\P^1)$ (with $\a$ as before) generate $\pi_2(D_u)$ and moreover, if $\a$ and $\b$ are such that $\a+\b$ is also a root, then $[\Sigma_{\a+\b}]=[\Sigma_\a]+[\Sigma_\b]$ in $\pi_2(D_u)$. Also, if there exists $\g\in\Phi\cap-\Phi$ such that $\a=\b+\g$ then $[\Sigma_\a]=[\Sigma_\b]$. Finally, we remark that the spheres $\Sigma_\a$ are superhorizontal for each $\a\in\Pi_\qg$. 
Then (see \cite{BR}),
\begin{prop}\label{pitwo}
$\pi_2(D_u)\simeq\Z^{k}$ and a basis of this $\Z$-module is given by $\{[\Sigma_\a]\,|\,\a\in\Pi_\qg\}$.
\end{prop}

Define $\beta_{i}:=\a_{r_0+\dots+r_i}\in\Pi_\qg$ and for any $i< j$, set $\beta_{i,j}:=\sum_{s=r_0+\dots+r_i}^{r_0+\dots+r_j}\a_s$.
Consider now the projection map $\pi_u:D_u\fd X_u$. For any $\a\in\Pi_\qg$, its restriction ${\pi_u}_{|\Sigma_\a}$ has degree one. Moreover, $\pi_2(X_u)\simeq\Z$ and for the right choice of a generator $[T]$, $[\pi_u(\Sigma_{\b_{j}})]=(-1)^j[T]$ (in fact, ${\pi_u}_{|\Sigma_{\b_{j}}}$ is alternately holomorphic and antiholomorphic). As a consequence, the morphism ${\pi_u}_\star$ induced by $\pi_u$ on the second homotopy groups is given by
$$\begin{array}{rcl}
{\pi_u}_\star : \pi_2(D_u)\simeq \Z^k & \fd & \pi_2(X_u)\simeq\Z\\
\sum_{i=0}^{k-1} a_{\b_{i}}[\Sigma_{\b_{i}}] & \longmapsto & \sum_{i=0}^{k-1} (-1)^i a_{\b_{i}} [T]\;\;.
\end{array}
$$
The central fiber $F$ of $\pi_u:D_u\fd X_u$ is ${\rm P}({\rm U}(p)\times {\rm U}(q))/{\rm P}(\Pi_{i=0}^k{\rm U}(r_i))$ which is nothing but ${\PU}(p)/{\rm P}(\Pi_{i\ {\rm even}}{\rm U}(r_i))\,\times\, {\PU}(q)/{\rm P}(\Pi_{i\ {\rm odd}}{\rm U}(r_i))$. It is hence isomorphic to the product of two flag manifolds:
applying Proposition~\ref{pitwo}, we get $\pi_2(F)\simeq\Z^{k-1}$ and by the above combinatorial
description, it is easily seen that $\pi_2(F)\simeq{\rm Ker}\,{\pi_u}_\star$, a basis of whom is given by
$$\bigl\{[\Sigma_{\b_{i}}]+[\Sigma_{\b_{i+1}}]\,|\,i=0,\dots,k-2 \bigr\}=\bigl\{[\Sigma_{\b_{i,i+1}}]\,|\,i=0,\dots,k-2 \bigr\}\;\;.
$$
Remark that each representative of the latter basis is contained in the central fiber of $\pi_u$. 

Since $F$ is also the central fiber of the fibration $\pi$ of the period domain $D$ over the noncompact symmetric space $X$, it is a deformation retract of $D$ and we deduce the
\begin{coro}
$\pi_2(D)\simeq\pi_2(F)\simeq\Z^{k-1}$.
\end{coro}

\subsubsection{Proof of Theorem~\ref{horiz}}

Proposition~\ref{petit} and the results of Section~\ref{bur}, enable us
to state the following theorem, of which Theorem~\ref{horiz} is a
particular case:
\begin{theo}\label{hori}
Let $E=\oplus_{i=0}^k E^i$ be a $\C$-VHS and set $r_i=\rk E^i$. Let
$D$ be the $\PU(p,q)$ period domain associated to $E$. Then,  for each
$i\in\{0,\dots,k-2\}$ such that $r_{i+1}\geq 2$, there exists a
superhorizontal representative $\psi:S^2\fd D$ of
$[\Sigma_{\b_{{i},{i+1}}}]$. In particular, when
the hypotheses of Theorem~\ref{cartoun} are not met, $\pi_2(D)$ can be generated by superhorizontal 2-spheres. 
\end{theo}

\begin{demo}
If $r_{i+1}\geq 2$,  $\{\b_i,\b_{i,i+1},\b_{i+1}\}$ are the simple roots of a Lie subalgebra of $\sg\lg(\E)$ isomorphic to $\sg\lg(4,\C)$ which is the complexification of a subalgebra of $\gg_0$ isomorphic to $\sg\ug(2,2)$. As in the case of $\P^1$ in Section~\ref{bur}, this defines a totally geodesic embedding of the period domain ${\mathcal D}=\PU(2,2)/{\rm P}(\U(1)\times\U(2)\times\U(1))$ in $D$. Geometrically, if $(e^i_1,\ldots,e_{r_i}^i)$ is a basis of $\E^i$, and if ${\mathbb W}=\la e^{i}_{r_{i}}\ra\oplus\la e_1^{i+1},e_{r_{i+1}}^{i+1}\ra\oplus\la e^{i+2}_1\ra$, then $\PU({\mathbb W},h_{|{\mathbb W}})\subset \PU(\E,h)$ is the copy of $\PU(2,2)$ corresponding to this embedding. 

Identifying ${\mathcal D}$ with its image in $D$, the superhorizontal distribution of ${\mathcal D}$ satisfies $H_{{\mathcal D}}\subset {H_{D}}_{|{\mathcal D}}$ and we have $\Sigma_{\b_{i,i+1}}\subset {\mathcal D}$. Hence we can apply Proposition~\ref{petit} in ${\mathcal D}$.
\end{demo}

\begin{rema} For a general period domain $D=G_0/V_0$, we can use this method to try to produce superhorizontal representatives of some elements of $\pi_2(D)$. Let us rephrase the proof we just gave. The group $\pi_2(D)$ is generated by vertical spheres.  More precisely, using the notation of the end of section~\ref{perdo}, if $\gg^\a\subset\gg_2\subset\gg$ is a root-space and $\gg^\a\oplus\gg^{-\a}\oplus[\gg^\a,\gg^{-\a}]$ is the corresponding copy of $\sg\lg(2,\C)$ in $\gg$, then exponentiating its real form $\sg\ug(2)$ gives a 2-sphere $\Sigma_\a$ in $D$, and these spheres generate $\pi_2(D)$. Now, if there is a Lie algebra homomorphism $\sg\lg(4,\C)\fd \gg$ which maps the root space $\sg\lg^{e_1-e_4}$ of $\sg\lg(4,\C)$ to $\gg^\a$ and the root spaces $\sg\lg^{e_1-e_2}$, $\sg\lg^{e_1-e_3}$, $\sg\lg^{e_2-e_4}$ and $\sg\lg^{e_3-e_4}$ to $\gg_1\subset\gg$, then as before there will be a copy of the period domain $\PU(2,2)/{\rm P}(\U(1)\times\U(2)\times\U(1))$ in $D$, with the right superhorizontality property and containing $\Sigma_\a$. Therefore Proposition~\ref{petit} yields a superhorizontal representative of the generator corresponding to the root-space $\gg^\a\subset \gg_2$.   

It is difficult to give a general criterion for the existence of these special homomorphisms $\sg\lg(4,\C)\fd \gg$. We list a few examples (there are many others) where the method can be applied at least partially, for classical matrix period domains of the type considered by Griffiths. In this case, $G_0$ is either ${\rm Sp}(n,\R)$ or $\SO(p,q)$ with $p$ or $q$ even. 
\begin{enumerate}
 \item for $D=\SO(2p,q)/(\U(p)\times\SO(q))$ with $p\geq 2$ and $q\geq 2$, $\pi_2(D)\simeq \Z$ can be generated by a superhorizontal sphere;
\item for $D=\SO(4,4)/(\U(1)\times\U(2)\times\SO(2))$, $\pi_2(D)\simeq \Z^3$ can be generated by superhorizontal spheres;
\item for $D=\SO(2p,2q_1+q_2)/(\U(q_1)\times\U(p)\times\SO(q_2))$ with $p\geq 2$, $q_1\geq 1$ and $q_2\geq 3$, $\pi_2(D)\simeq \Z^2$ can be generated by superhorizontal spheres;  
\item for $D={\rm Sp}(n,\R)/(\U(k)\times\U(n-2k)\times\U(k))$ with $k\geq 1$ and $n-2k\geq 2$, $\pi_2(D)=\Z^2$ and (at least) one generator can be taken to be a superhorizontal sphere;
\item  $D={\rm Sp}(n,\R)/(\U(k)\times\U(n-k))$ with $k\geq 1$ and $n-k\geq 3$, $\pi_2(D)\simeq \Z$ can be generated by a superhorizontal sphere.
\end{enumerate}
\end{rema}

\begin{rema} Concerning the possibility of generating the second homotopy groups of $\PU(p,q)$ period domains by horizontal spheres, there is one question which remains unanswered. In the notation of Theorem~\ref{hori}, we do not know whether $[\Sigma_{\b_{i,i+1}}]$ can be represented by a superhorizontal 2-sphere when $r_{i+1}=1$. 

To take the simplest example, it would be very interesting to know what happens for the period domain $D={\rm PU}(2,1)/{\rm P}(\U(1)\times\U(1)\times \U(1))$, whose second homotopy group is isomorphic to $\Z$. As pointed out in Remark~\ref{n=1}, our implementation of the $h$-principle in ${\rm PU}(2,n)/{\rm P}(\U(1)\times \U(n)\times \U(1))$ fails when $n=1$. 

Nevertheless, the arguments of the proof of Theorem~\ref{cartoun} imply that for any period domain $D$ of the form ${\rm PU}(p+q,1)/{\rm P}(\U(p)\times \U(1)\times \U(q))$ (with $p\geq 1$, $q\geq 1$), and for any holomorphic superhorizontal map $f:\tilde M\fd D$, we have $f_\star(\pi_2(M))=0$. 

Indeed, let us consider the holomorphic projection $\nu:D\fd D':={\rm PU}(p+q,1)/{\rm P}(\U(p)\times \U(q,1))$ and the projection onto the symmetric space $\pi:D\fd X:={\rm PU}(p+q,1)/{\rm P}(\U(p+q)\times \U(1))$. Then, for any 2-sphere $\sigma:S^2\fd \tilde M$, $\nu\circ f\circ\sigma:S^2\fd D'$ and $\pi\circ f\circ\sigma:S^2\fd X$ are homotopically trivial (the former by the arguments of Section~\ref{proofofcartoun}, the latter because $X$ is contractible).
Let $\tau_1:S^2\times [0,1]\fd D'$ be a homotopy between $\nu\circ f\circ\sigma$ and a constant map and $\tau_2:S^2\times [0,1]\fd X$ be a homotopy between $\pi\circ f\circ\sigma$ and a constant map. For any $(s,t)\in S^2\times [0,1]$, the flag $\{0\subset\tau_1(s,t)\subset \tau_1(s,t)\oplus\tau_2(s,t)^\perp\subset\C^{p+q+1}\}$ is an element of $D$ (orthogonality is understood with respect to the indefinite Hermitian form of signature $(p+q,1)$ on $\C^{p+q+1}$), although $\tau_1(s,t)$ and $\tau_2(s,t)^\perp$ are not necessarily orthogonal. In this way we get a homotopy between $f\circ\sigma$ and a constant map.
\end{rema}

\end{document}